%% file: 0---laplace-beltrami-2.7.tex
\documentclass[11pt]{article}
\usepackage{amssymb,amsmath,mathrsfs,amsfonts,graphicx,epsf} 
\usepackage{amsmath,amsfonts,amssymb,mathrsfs,amsthm,color}
\usepackage[english,activeacute]{babel}

\usepackage[latin1]{inputenc}
\usepackage[T1]{fontenc}
 \input xy
 \xyoption{all}

 \usepackage{esint}

 \usepackage{fullpage}

  \newcommand{\ffoot}[1]{}

\newcommand{\eucl}{e}

\newcommand{\bD}{\blacktriangle}

\renewcommand{\th}{\theta}

\newcommand{\f}{\mathfrak{f}}

\usepackage{mathrsfs}
\usepackage{graphicx}
\usepackage{amssymb}
\usepackage{color}
\newtheorem{theorem}{Theorem}
\newtheorem{corollary}{Corollary}
\newtheorem{lemma}{Lemma}
\newtheorem{proposition}{Proposition}
\newtheorem{definition}{Definition}
\newtheorem{remark}{Remark}

\newcommand{\bt}{\begin{theorem}}
\newcommand{\et}{\end{theorem}}
\newcommand{\bl}{\begin{lemma}}
\newcommand{\el}{\end{lemma}}
\newcommand{\bp}{\begin{proposition}}
\newcommand{\ep}{\end{proposition}}
\newcommand{\bc}{\begin{corollary}}
\newcommand{\ec}{\end{corollary}}
\newcommand{\bdeff}{\begin{definition}}
\newcommand{\edeff}{\end{definition}}
\newcommand{\brem}{\begin{remark}}
\newcommand{\erem}{\end{remark}}
\renewcommand{\r}[1]{(\ref{#1})}
\newcommand{\con}{{\mathcal C}}

\newcommand{\bi}{\begin{itemize}}
\newcommand{\iii}{\item}
\newcommand{\ei}{\end{itemize}}
\newcommand{\bd}{\begin{description}}
\newcommand{\ed}{\end{description}}

\newcommand{\bqn}{\begin{eqnarray}}
\newcommand{\eqn}{\end{eqnarray}}
\newcommand{\eqnn}{\nonumber\end{eqnarray}}

\newcommand{\nn}{\nonumber}
\newcommand{\ba}[1]{\begin{array}{#1}}
\newcommand{\ea}{\end{array}}

\newcommand{\R}{\mathbb{R}}
\newcommand{\Z}{\mathbb{Z}}
\newcommand{\N}{\mathbb{N}}

\newcommand{\lettre}{\beta}

\newcommand{\lam}{\lambda}
\newcommand{\g}{\gamma}
\newcommand{\al}{\alpha}
\newcommand{\eps}{\varepsilon}

\newcommand{\VecM}{\mathrm{Vec}(M)}

\newcommand{\Gq}{{\gg}_q}

\newcommand{\Zz}{\mathcal{Z}}
\renewcommand{\gg}{{\bf G}}

\newcommand{\HH}{{\bf (H0)}}
\def\R{{\mathbf{R}}}

\newcommand{\frp}[2]{\frac{\partial #1}{\partial #2}}

\def\F{{\mathbf{F}}}

\def\Z{{\mathbf{Z}}}

\newcommand{\grad}{\nabla}
\renewcommand{\div}{{\mbox{div}}}

\newcommand\bna{\begin{eqnarray*}}%équation non numérotée
\newcommand\ena{\end{eqnarray*}}
\newcommand\bnan{\begin{eqnarray}}%équation numérotée
\newcommand\enan{\end{eqnarray}}

\newcommand\bneq{\begin{eqnarray*}\left\lbrace \begin{array}{rcl}}%systeme equation sans numerotation
\newcommand\eneq{\end{array} \right.\end{eqnarray*}}
\newcommand\bneqn{\begin{eqnarray}\left\lbrace \begin{array}{rcl}}%systeme equation avec numerotation
\newcommand\eneqn{\end{array} \right.\end{eqnarray}}

\newcommand\bnp{\begin{proof}}%équation numérotée
\newcommand\enp{\end{proof}}
\newcommand\nor[2]{\left\|#1\right\|_{#2}}%la norme
\newcommand\norL[1]{\left\|#1\right\|_{L^2}}%la norme L^2
\newcommand\sgn{\textnormal{sgn}}
\newcommand\Tu{\mathbb{T}}

\newcommand\grando[1]{\mathcal{O}(#1)}

\newcommand\limvar[2]{\underset{#1\to #2}{\lim}}
\newcommand{\comp}{W}

\newcommand{\ugo}[1]{#1}%{\color{blue}{#1}}}

\begin{document}

\title{{\bf The Laplace-Beltrami operator in almost-Riemannian Geometry}}

\author{ 
 Ugo~Boscain, \\[2mm]
CNRS,  Centre de Math\'ematiques Appliqu\'ees, \'Ecole Polytechnique\\[2mm]
Route de Saclay, 91128 Palaiseau Cedex, France,\\[2mm]
and equipe INRIA GECO\\[2mm]
        {\tt ugo.boscain@cmap.polytechnique.fr}\\[4mm]
  Camille Laurent\\[2mm]
CNRS,  Centre de Math\'ematiques Appliqu\'ees, \'Ecole Polytechnique\\[2mm]
Route de Saclay, 91128 Palaiseau Cedex, France,\\[2mm]
        {\tt camille.laurent@cmap.polytechnique.fr}
 }

\maketitle

\begin{abstract}
{   Two-dimensional almost-Riemannian structures are generalized Riemannian structures on surfaces for which a local orthonormal frame is given by a Lie bracket generating pair of vector fields that can become collinear. 
Generically, the singular set is an embedded one dimensional manifold and there are three type of points: Riemannian points 
where the two vector fields are linearly independent, Grushin points where the two vector fields 
are collinear but their Lie bracket is not and tangency points where the two vector fields and 
their Lie bracket are collinear and the missing direction is obtained with one more bracket. Generically tangency points are isolated.

In this paper we  study the Laplace-Beltrami operator on such a structure. In the case of a compact orientable surface without tangency points, we prove that the Laplace-Beltrami operator is essentially self-adjoint and has discrete spectrum.  As a consequence a quantum particle in such a structure cannot cross the singular set and the heat cannot flow through the singularity. This is an interesting phenomenon since when approaching  the  singular set (i.e. where the vector fields become collinear),  all Riemannian quantities \ugo{explode}, but geodesics are still well defined and can cross the singular set without singularities. 

This phenomenon appears also in sub-Riemannian structure which are not equiregular i.e. in which the grow vector depends on the point. We show this fact by analyzing the Martinet case.\footnote{\ugo{This research has been supported by the European Research Council, ERC StG 2009 {GeCoMethods}, contract number 239748, by the ANR Project GCM, program {Blanche}, project number NT09-504490 and by the DIGITEO project CONGEO.}}

}

\end{abstract}

\section{Introduction}

A $2$-dimensional Almost Riemannian  Structure ($2$-ARS for short) is a generalized Riemannian structure  that can be defined locally by a pair of smooth vector fields on a $2$-dimensional manifold $M$, satisfying the H\"ormander condition. These vector fields play the role of an orthonormal frame. 

Let us denote by $\bD(q)$  the linear span of the two vector fields at a  point $q$.  Where  $\bD(q)$ is $2$-dimensional, the corresponding metric is Riemannian. Where $\bD(q)$ is $1$-dimensional, the corresponding Riemannian metric is not well defined, but thanks to the H\"ormander condition one can still define the Carnot-Caratheodory distance between two points, which happens to be finite and continuous. 

$2$-ARSs were introduced in the context of hypoelliptic operators \cite{FL1,grushin1}, they appeared in  problems of  population transfer in quantum systems  \cite{q4,BCha,q1}, and have applications to  orbital transfer in space mechanics  \cite{BCa,tannaka}. 
$2$-ARSs  are a particular case of rank-varying sub-Riemannian structures (see for instance \cite{bellaiche,jean1,jean2}). 

\ugo{Generically (i.e for an open and dense subset of the set of all 2-ARSs, in a suitable topology), the singular set  $\Zz$, where $\bD(q)$ has dimension $1$,  is a $1$-dimensional embedded submanifold and there are three types of points: Riemannian points, Grushin points where  $\bD(q)$  is $1$-dimensional and dim$(\bD(q)+[\bD,\bD](q))=2$ and  tangency points where dim$(\bD(q)+[\bD,\bD](q))=1$ and  the missing direction is obtained with one more bracket. 
One can easily show that  at  Grushin points $\bD(q)$  is transversal to $\Zz$. Generically, at tangency points $\bD(q)$ is tangent to $\Zz$ and tangency points are isolated.} Normal forms at Riemannian, Grushin and tangency points were established in \cite{ABS} and are described in Figure 1. %\ref{fig-prenormal}.

$2$-ARSs present very interesting phenomena.  For instance, geodesics  can pass through the singular set, with no singularities even if all Riemannian quantities (as for instance the metric, the Riemannian area, the curvature) explode while approaching $\Zz$. Moreover  the presence of a singular set permits the conjugate locus  to be nonempty even if the Gaussian curvature is always negative, where it is defined (see \cite{ABS}).
See also  \cite{ABS,euler,BCGS,high-order} for  Gauss--Bonnet-type formulas and for a classification of 2-ARSs from the point of view of Lipschitz equivalence.

%However, tangency points are far to be deeply understood.  An open question is the convergence or the divergence of the integral of the geodesic curvature on the boundary of  a tubular  neighborhood  of the singular set, close to a tangency point. This question arose in the proof of the Gauss--Bonnet theorem given in \cite{euler}. In that paper, thanks to numerical simulations, the authors conjecture  the divergence of such integral.
 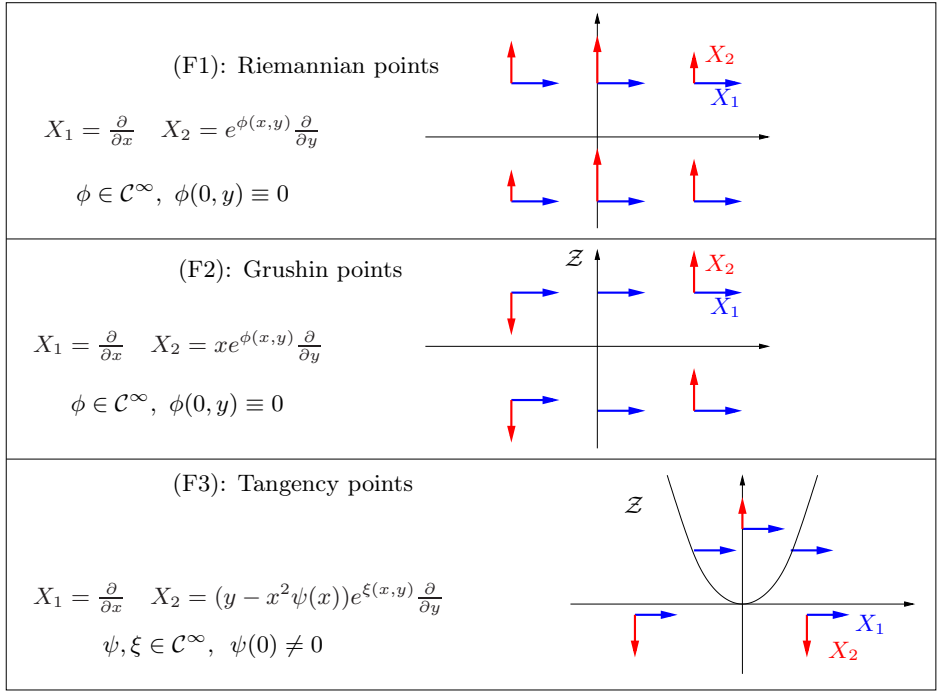
\begin{figure}[h!]
\begin{center}
\input{fig-intro-ugo.pstex_t}
\caption{The local representations established in \cite{ABS}.  Tangency points are the most difficult to handle due to the fact that the asymptotic of the distance is different from the two sides of the singular set (see \cite{BCGJ,kresta}). }
\label{fig-prenormal}
\end{center}
\end{figure}

In this paper we study the Laplace-Beltrami operator in a $2$-ARS. The main point is that the first order terms explode as a consequence of the explosion of the area on the singular set $\Zz$. In particular we are interested in its self-adjointness.  This is crucial to understand the evolution of a quantum particle and of the heat flow in a 2-ARS.

On $M\setminus\Zz$, this operator is defined as the divergence of the gradient.  The almost-Riemannian gradient can be defined with no difficulty since it involves only the inverse metric which is well defined  (by continuity) even on the singular set.   
More precisely, on a  Riemannian manifold $(M,g)$, $\grad$ is the unique operator from  $\mathcal{C}^\infty(M)$ to Vec$(M)$ satisfying $g_q(\grad~\phi(q),v)=d\phi_q (v)$, $\forall~q\in M,~v\in T_qM$. In the two-dimensional case, if in an open set $\Omega\subset M$ an orthonormal frame for $g$ is given by two vector fields $X_1$ and $X_2$,  then we have on  $\Omega$,
\bqn
\grad (\phi)= X_1(\phi) X_1+X_2(\phi) X_2,
\label{ars-grad}
\eqn
where by $X_i(\phi)$ we mean the Lie derivative of $\phi$ in the direction of $X_i$ ($i=1,2$). This last formula can be used to define the gradient even where $X_1$ and $X_1$ are not linearly independent. In this paper, we use formula \r{ars-grad} to define the gradient of a function on a 2-ARS.

However the divergence of a vector field requires a notion of area and  in a 2-ARS the natural area (i.e. the Riemannian one) explodes while $q$ approaches the singular set.  
Let us recall that the divergence of a vector field $X$ on a Riemannian manifold $(M,g)$ is the unique function satisfying   
$\div X d\omega=L_X d\omega$ where $d\omega$ is the Riemannian volume form. In coordinates $\div X=\frac{1}{\omega}\sum_{i,j}\partial_i\left(\omega X^j\right)$, 
where $\omega=\sqrt{\text{det} g}$. 

 To see how the Riemannian quantities (and in particular $d\omega$) explode while approaching  the singular set $\Zz$, choose a local orthonormal frame for the 2-ARS of the form  $X_1=(1,0)$ and $X_2=(0,f(x,y))$,  where $f$ is some smooth function. This is always possible, see for instance \cite{ABS}.   On $M\setminus\Zz$,  for the metric $g$, the area element $d\omega$, the curvature $K$ and the gradient of a smooth function one has:
\bqn
g(x,y)&=&\left(\ba{cc}1&0\\0&\frac{1}{f(x,y)^2}\ea \right),\nn\\
d\omega(x,y)&=&\frac{1}{|f(x,y)|} dx\, dy,\nn\\
K(x,y)&=&  \frac{f(x,y) \partial_{x}^2  f(x,y)-2
   \partial_{x}f(x,y)^2}{f(x,y)^2},    \nn\\
\grad (\phi)&=& \partial_{x}(\phi)\partial_{x}+f^2 \partial_{y}(\phi)\partial_{y}.\nn
\eqn
and for a vector field $Y=(Y^1,Y^2)$ one has
\bqn
\div(Y)&=&|f|\left(     \partial_{x}\left(\frac{1}{|f|} Y^1\right) + \partial_{y}\left(\frac{1}{|f|} Y^2\right)\right)\nn\\
&=&\partial_{x}Y^1+\partial_{y}Y^2-
\frac{\partial_{x}f}{f}Y^1
-\frac{\partial_{y}f}{f} Y^2
\eqn
which is not well defined on $\Zz$ even on vector  fields which are linear combinations of $X_1$ and $X_2$.    In particular it is not well defined on $X_1$. As a consequence the Laplace Beltrami operator presents some singularities in the first order terms. (The choice of the area does not affect the principal symbol of the operator). More precisely, one has
\bqn
\Delta(\phi)=\partial_{x}^2\phi+f^2 \partial_{y}^2\phi -
\frac{\partial_{x}f}{f}\partial_{x}\phi
+f(\partial_{y}f) (\partial_{y}\phi).
\eqn

The simplest case is the well known Grushin metric, which is the 2-ARS on $\R^2$ for which an orthonormal basis is given by, 
\bqn
X_1=
\left(\ba{c}
1\\0
\ea
\right),~~~X_2=
\left(\ba{c}
0\\x
\ea
\right).
\eqn
Here the singular \ugo{set $\Zz$} is the $y$ axis and on $\R^2\setminus\Zz$ the Riemannian metric, the Riemannian area and the Gaussian  curvature are given respectively by:
\bqn
g=\left(\ba{cc}
1&0\\0&\frac1{x^2}
\ea
\right),~~~
d\omega=\frac{1}{|x|}dx\, dy,~~~K=-\frac{2}{x^2}.
\eqn
It follows that 
%given a smooth function $\phi: \R^2\to\R$ and for a smooth vector field   $X: \R^2\to\R^2$, we have 
%{\tt Far vedere che le geodetiche sono ben definite: due ficure: partendo dal set e da prima} 
%\bqn
%\grad(\phi)=......,\div X=......., \eqn
%and 
the Laplace Beltrami operator is given by
\bqn
\Delta \phi:=\div(\grad(\phi))=\big(\partial_{x}^2+x^2\partial^2_{y}-\frac{1}{x}\partial_{x}\big)\phi,
\label{eq-LBgrushin}
\eqn
which is the basic model of Laplace Beltrami operator on 2-ARSs and whose self-adjointness will be studied in this paper. 
Notice that the operator "sum of squares",
\bqn
\bar\Delta=X_1^2+X_2^2=\partial_{x}^2+x^2\partial_{y}^2
\eqn
has been deeply studied in the literature. This operator (which in this case is the principal symbol of $\Delta$) can be obtained as the divergence of the gradient by taking as area, the standard Euclidean area on $\R^2$.
However this operator does not give any information about $\Delta$ because of the diverging first order term. \ugo{Moreover in 2-ARS 
an intrinsic Laplacian with non-diverging terms cannot be  defined naturally, since an intrinsic area which does not diverge on the singular set  is not known.} Here by intrinsic we mean an area which depends only on the 2-ARS and not on the choice of the coordinates and of the orthonormal frame.

Notice that even if all Riemannian quantities are not defined on $M\setminus\Zz$, classical geodesics do. See Section \ref{ss-geo-grushin} for the explicit expressions of geodesics in the Grushin plane.
 Another interesting feature of the Grushin plane, is the fact that a bounded open set intersecting $\Zz$  has finite diameter but infinite area.

The main purpose of this paper is to study the following question:

\bd
\iii[{[Q]}] Let $M$ be a 2-D manifold endowed with a 2-ARS, and let $\Delta$ be the corresponding Laplace-Beltrami defined on $C^\infty_0(M\setminus\Zz)$. Is $\Delta$ essentially self-adjoint?
\ed
Notice that  a priory one expects a negative answer to this question, since, as explained below, a positive answer would imply that 
neither \ugo{the heat flow, neither a quantum particle can pass through $\Zz$, while classical geodesics cross it with no singularities.}

Our main result is an unexpected positive answer to the question {\bf [Q]} in the case in which there are no tangency points and in the case in which $M$ is compact. More precisely we have the following.
\bt
\label{t-1}
Let $M$ be a 2-D compact  orientable manifold endowed with a  2-ARS.  Assume that 
\bd
\iii[{[HA]}] the singular set $\Zz$ is an
embedded one-dimensional 
submanifold of
$M$;
\iii[{[HB]}] for every $q\in M$, $\bD(q)+[\bD,\bD](q)=T_qM$.
\ed
 Let $d\omega$ be the corresponding Riemannian area and $\Delta$ the corresponding Laplace-Beltrami operator, both  defined on $M\setminus\Zz$.  Then  we have the following facts
\begin{enumerate}
\item $\Delta$ with domain $C^\infty_0(M\setminus\Zz)$  is essentially self-adjoint on $L^2(M,d\omega)$. 
 \item 
The domain of $\Delta$ (the self-adjoint extension denoted by the same letter) is given by
\bna
D(\Delta)=\left\lbrace u\in L^2(M,d\omega) \left| \Delta_{D,g} u \in L^2(M,d\omega)\right.\right\rbrace 
\ena
 where $\Delta_{D,g} u$ is $\Delta u$ seen as a distribution in $\Omega=M\setminus \Zz$. 
\item The resolvent $(-\Delta+1)^{-1}$ is compact and therefore its spectrum is discrete and consists of eigenvalues with finite multiplicity.
\end{enumerate}
\et
Hypothesis {\bf HA} is generic (see \cite{ABS}). \ugo{Hypothesis}  {\bf HB} implies that every point is either a  Riemannian  point or a Grushin point  (see  Figure 1).  %\ref{fig-prenormal}) . 
The hypotheses that $M$ is compact and that there no tangency points are technical and the same result should hold in much more general situations.
The orientability of the manifold can be weakened. Indeed it is only necessary that each connected component of $\Zz$ (which is diffeomorfic to $S^1$) admit an open  tubular neighborhood diffeomorfic to $(-1,1)\times S^1$. See Proposition \ref{p-globalizziagite} and \ref{r-or}.

\brem
Notice that 2-ARSs satisfying {\bf HA} and {\bf HB} do exist. See for instance \cite{ABS} where such a 2-ARS has been built on compact orientable manifolds of any genus. Notice moreover that a 2-ARS  satisfying {\bf HA} and {\bf HB}
are structurally stable in the sense that small perturbations of the local orthonormal frames in the $C^\infty$ norm, do not destroy conditions {\bf HA} and {\bf HB}.
\erem
Theorem \ref{t-1} has a certain number of implications. 
First it implies that a quantum particle, \ugo{localized for $t=0$} on a connected component of $M\setminus\Zz$, {   remains localized in this connected component for any time $t\in \R$. The same phenomenon holds for the wave or the heat equation. Indeed, the essential self-adjointness means that our operator can be naturally and uniquely extended (by taking its closure) to a self-adjoint operator without adding any additional boundary condition. But, in our case, one possible extension is to take an extension for the operator defined on each connected component (for instance the Friedrichs extension) and to ``concatenate'' them. This possible extension is the one which separates the dynamics. For instance this is what we obtain when  one defines  the Friedrichs extension of the usual Laplace operator on $\R$ defined on function in $C^{\infty}_0(\R \setminus\{0\})$. We obtain two distinct dynamics with Dirichlet boundary condition at $0$. Yet, this is of course not the only possible dynamic because the operator defined on function $C^{\infty}_0(\R \setminus\{0\})$ is not essentially self-adjoint. So, the essential self-adjointness means that the unique self-adjoint extension of our operator is the one that separates the dynamics on each connected components. That means that it is not necessary to add any boundary condition: the explosion of the area naturally acts as a barrier which prevents the crossing of the degeneracy zone by the particules. }

More precisely we have the following.
\bc
With the notations of Theorem \ref{t-1}, consider the unique solution $u$ of the Schrödinger equation (according to the self-adjoint extension defined in the previous theorem),
\bneqn
i\partial_t u+\Delta u&=&0\\
u(0)&=&u_0\in  L^2(M,d\omega)
\eneqn
with $u_0$ supported in a connected component $\Omega$ of $M\setminus\Zz$. Then, $u(t)$ is supported in $\Omega$  for any $t\geq 0$.
The same holds for the solution of the heat or for the solution of the wave equation. 
\ec

Second it is well known that in Riemannian geometry one can relate  properties of the Riemannian distance to those of the corresponding heat kernels. For instance we have that 
{   
\bt[Varadhan, Neel, Stroock] Let $M$ be a compact, connected, smooth Riemannian manifold and $d$ the corresponding Riemannian distance.  Let $p_t(q_1,q_2)$ be the heat kernel of the heat equation  $\partial_t \phi=\frac12\Delta \phi$, where $\Delta$  is the Laplace-Beltrami operator. Define 
$$
E_t(q_1, q_2):=-t \log p_t (q_1, q_2)
$$
We have:
\bd
\iii[1)] $E_t(q_1, q_2)\to \frac12 d(q_1,q_2)^2$, uniformly on $M\times M$ as $t\to 0$, 
\iii[2)] Let Cut$(q_1)\subset M$ the cut locus from $q_1$. Then
$q_2\not\in$Cut$(q_1)$ if and only if $\lim_{t\to 0}\nabla^2 E_t(q_1,q_2)=\frac12 \nabla^2 d(q_1,q_2)^2$, while $q_2\in$Cut$(q_1)$ if and only if $\limsup_{t\to 0}\|\nabla^2 E_t(q_1,q_2)\|=\infty$, where $\|.\|$ is the operator norm.
\ed
\et
The first result is due to Varadhan \cite{varadhan}, the second one to Neel and Stroock \cite{neel,neel-stroock}.
}
%%%%%%%%%%%%%%%

%%%%%%%%%%%%%%%
Both these results do not depend upon the fact that one is using the Laplace-Beltrami operator or the principal symbol of the operator in a fixed frame to construct the heat kernel.

The situation is very different in 2-ARG. Both these results are false for the Laplace Beltrami operator. Indeed the distance between two points belonging to two different connected components of $M\setminus \Zz $ is finite while if $q$ belongs to a connected component $W$ of $M\setminus \Zz $  then $p_t(q,\cdot)$ is supported in $W$. 

However a result in the spirit of the one of Varadhan has been obtained by Leandre in \cite{leandre} for the operator ``sum of squares''.  Hence the Laplace-Beltrami operator defined above, has quite different properties with respect to the operator sum of squares. 
The last one is not intrinsic, but however keeps tracks of intrinsic quantities as the almost-Riemannian distance. In particular the corresponding heat flow crosses the set $\Zz$ which is not the case for the Laplace-Beltrami operator.  
For other relations among the heat kernel and the distance in sub-Riemannian geometry see \cite{davide-solo}.

  {To prove Theorem \ref{t-1} we start by analyzing the Grushin case \r{eq-LBgrushin}. We first compactify in the $y$ variable by considering it on $\R_x\times \Tu_y$.

By setting $f=\sqrt{|x|}g$ we are reduced to study the essential self adjointness of the following operator considered on $L^2$ with the usual euclidian metric
\bqn
L=\partial_{x}^2+x^2\partial_{y}^2-\frac{3}{4}\frac{1}{x^2}.
\eqn 
 If for a moment we forget the term $x^2\partial_{y}^2$ which is not relevant for the selfadjointness of this operator (and becomes a non-positive potential $-k^2x^2$ after performing Fourier transform in $y$), we are reduced to study the operator $\partial_{x}^2-\frac{3}{4}\frac{1}{x^2}$, which is well known in the literature. Indeed we have the following result 
\bp
The operator $-\partial_x^2 +\frac{c}{|x|^2}$ defined on $L^2(]0,+\infty[)$ with domain $C^{\infty}_0(]0,+\infty[)$ is essentially self-adjoint if and only $c\geq \frac{3}{4}$.
\ep
There are different proof of this result in the literature (see \cite{ReedSimon} Chapter X for an introduction), often leading to some stronger statement (for instance for potential $V\geq \frac{3}{4x^2})$) and generalizations to higher dimensions. 

The rest of the proof for an almost Riemannian structure \ugo{consists in generalizing this result for a normal form around a connected component of the singular set}. The main tools are Kato inequality and perturbation theory.

\brem
Notice  that, as remarked in the introduction of \cite{PGGrellier},  solutions of the Schr\"odinger equation,  with the Laplacian defined as sum of square $\partial_x^2+x^2 \partial_y^2$, display a total lack of dispersion. It would be interesting to understand if the same holds for the Laplace-Beltrami operator defined by (\ref{eq-LBgrushin}). It would be also interesting to study the behaviour of a wave packet moving towards  the singular set, to understand if there is a reflection or a dispersion. Since in Theorem \ref{t-1} we prove that the spectrum is discrete, it seems natural to expect a reflection.
\erem

The structure of the paper is the following. In section~\ref{s-basdef} we   briefly recall the notion of almost-Riemannian structure. An original result is given in Proposition \ref{p-globalizziagite},  in which we globalize the normal form of type $(F_2)$  around a connected component of $\Zz$ in the compact orientable case and under the assumption that there are no tangency points. In Section \ref{s-grushin} we analyze the Grushin case, both for what concerns geodesics and the self-adjointness of the Laplace-Beltrami operator.
The main result, is proven in Section~\ref{s-camillo}.

 As a byproduct of our studies, we obtain that the Laplace-Beltrami operator is not essentially self-adjoint in the case in which an orthonormal basis is given by the vectors $(1,0)$, $(0,|x|^\al)$ with $\al\in(0,1)$ (see remark \ref{r-alfa}). This fact motivates some remarks on metrics of this kind (which do not enter in the standard framework of almost-Riemannian geometry, since the vector fields are not smooth). For instance we show that  in this case, even if  the area element explodes, the area of a  bounded open set intersecting $\Zz$  is finite (which is false for 2-ARSs). Moreover we show that there exist regular curves passing through $\Zz$ with a velocity not belonging to the span of $X_1$ and $X_2$, having finite length. 
See Section \ref{s-sqrt}.

 Finally, in the appendix, we discuss the Martinet case to show that this kind of phenomenon 
appears also in sub-Riemannian structure of constant rank but which are not equiregular i.e. in which the grow vector depends on the point (see for instance \cite{Montgomery}). In this case the  role of the Laplace-Beltrami operator is played by the intrinsic sub-Riemannian Laplacian defined via the Popp measure (see \cite{HeatKernel,Montgomery}).  
While aproaching the Martinet surface  the  Popp volume explodes and  the first order coefficient of the  sub-Riemannian Laplacian does as well. 
%In this Theorem the variable $x$ and $z$ have been compactified by simplicity

\bt
\label{t-martinet}
Consider the sub-Riemannian structure in $M=\Tu_x\times\R_y\times\Tu_z$ for which an orthonormal basis is given by 
$X_1=(
1,0,\frac{y^2}{2}),$ $X_2=
(0,1,0)$. Then the corresponding intrinsic sub-Riemannian Laplacian which is given by $\Delta_{sr}=(X_1)^2+(X_2)^2-\frac{1}{y}X_2$ with domain $C^\infty_0(M\setminus\{y=0\})$  is essentially self-adjoint on $L^2(M,d\omega)$, where $d\omega=\frac{1}{|y|}dx\,dy\,dz$ is the Popp volume.
\et

\section{Basic Definitions}
\label{s-basdef}
In this section we recall some basic definitions  in the framework
 of  2-ARSs following  \cite{ABS,euler}.

Let $M$ be a  smooth  surface without boundary.  Throughout the paper, unless specified,  manifolds are smooth (i.e., $C^{\infty}$) and without boundary; vector fields  and differential forms are smooth. The set of smooth vector fields
 on $M$ is denoted by $\VecM$.  The circle is denoted by $S^1$ or $\Tu$ depending on the context. $C^{\infty}_0$ denotes the set of smooth functions with compact support.
 
 \subsection{2-Almost-Riemannian Structures}
 \begin{definition}
 \label{d-ARS}
 A {\it $2$-dimensional almost-Riemannian structure} (2-ARS) is a triple
${\mathcal S}=
(E, \f,\langle\cdot,\cdot\rangle)$
 where 
  $E$ is a  vector bundle of rank $2$ over $M$ and $\langle\cdot,\cdot\rangle$ is a Euclidean structure on $E$, that is, $\langle\cdot,\cdot\rangle_q$ is a   scalar product 
on $E_q$ smoothly depending on $q$.  Finally
   $\f:E\rightarrow TM$ is a morphism of vector bundles,    
i.e., {\bf (i)}  the diagram 
$$
\xymatrix{
 E  \ar[r]^{\f} \ar[dr]_{\pi_E}   & TM \ar[d]^{\pi}            \\
 & M                          
}    
$$
 commutes, where  $\pi:TM\rightarrow M$ and $\pi_E:E\rightarrow M$ denote
  the canonical projections and {\bf (ii)} $\f$ is linear on fibers.

Denoting by $\Gamma(E)$ the $C^\infty(M)$-module of 
smooth sections on $E$, and by 
 $\f_*:\Gamma(E)\rightarrow \VecM$, the map $\sigma\mapsto\f_*(\sigma):=\f\circ\sigma$. We require that  the submodule of Vec$(M)$ given by $\bD=\f_*(\Gamma(E))$ to be bracket generating, i.e.,
  $Lie_q(\bD)= T_qM$ for every $q\in M$. 
\end{definition}
Here Lie$(\bD)$ is  the smallest Lie subalgebra of Vec(M) containing $\bD$ and 
Lie$_q(\bD)$ is the linear subspace of $T_qM$ whose elements are evaluation at $q$ of elements belonging to Lie$(\bD)$.
The condition that $\bD$ satisfies the Lie bracket generating assumption is known also as the \ugo{H\"ormander} condition.

{   We say that a 2-ARS $(E, \f,\langle\cdot,\cdot\rangle)$ is orientable if E is orientable as a vector bundle. Notice that one can build non-orientable 2-ARSs on orientable manifolds and orientable 2-ARSs on non-orientable manifolds. See \cite{ABS} for some examples.
We say that a 2-ARS $(E, \f,\langle\cdot,\cdot\rangle)$  is trivializable if $E$ is isomorphic to the trivial bundle $M 
\times \R^2$. A particular case of 2-ARSs is given by Riemannian surfaces. In this case $E=TM$ and $\f$ is the identity.
}

Let ${\mathcal S}=(E,\f,\langle\cdot,\cdot\rangle)$ be a  2-ARS on a surface $M$.
We  denote by $\bD(q)$  the linear subspace $\{V(q)\mid  V\in \bD\}=\f(E_q)\subseteq  T_q M$. 
The set of points in $M$ such that $\dim(\bD(q))<2$ is called {\it singular set} and denoted by $\Zz $.
Since $\bD$ is bracket generating, the subspace $\bD(q)$ is nontrivial for every $q$ and $\Zz$ coincides with the set of points $q$ where $\bD$ is one-dimensional.

 The Euclidean structure on $E$ allows to define a symmetric positive definite $C^\infty(M)$-bilinear form on the submodule $\bD$ by
 \bqn 
 G:\bD\times\bD\to C^\infty(M)\nn\\
 G(V,W)=\langle\sigma_V,\sigma_W\rangle\nn
 \eqn
 where $\sigma_V,\sigma_W$ are the unique\footnote{the uniqueness is consequence of the fact that we assume $E$ of rank two and $\bD$ Lie bracket generating.}
 sections of $E$ satisfying $\f\circ\sigma_V=V,\f\circ\sigma_W=W$. 
 
 At points $q\in M$ where $\f|_{E_q}$ is an isomorphism, i.e. on $M\setminus\Zz$,  $G$ is a tensor and the value $G(V,W)|_q$ depends only on $V(q), W(q)$. In this case $G$ defines a Riemannian metric $g$ via
 $$
 g_q(v,w)=G(V,W)(q),~~v,w\in T_q M, 
 $$
 where $V$ and $W$ are two vector fields such that $v=V(q)$ and $w=W(q)$.
 
 This is no longer true at points  $q$ where $\f|_{E_q}$ is not injective.

If $(\sigma_1,\sigma_2)$ is an orthonormal frame for $\langle\cdot,\cdot\rangle$ on
 an open subset $\Omega$ of $M$, an {\it  orthonormal frame for $G$} on $\Omega$ is 
 given by  $(\f\circ\sigma_1,\f\circ\sigma_2)$.  Orthonormal frames are systems of local generators of $\bD$.

For every $q\in M$  
and every $v\in\bD(q)$ define
$$
\Gq(v)=\inf\{\langle u, u\rangle_q \mid u\in E_q,\f(u)=v\}.
$$
An  absolutely continuous curve $\g:[0,T]\to M$  is  admissible for ${\mathcal S}$ 
if   there exists a measurable essentially bounded function 
$$[0,T]\ni t\mapsto u(t)\in E_{\g(t)},
$$ called {\it control function} such that 
$\dot \g(t)=\f(u(t))$  for almost every $t\in[0,T]$. 
Given an admissible 
curve $\g:[0,T]\to M$, the {\it length of $\g$} is  
\bqn
\ell(\g)= \int_{0}^{T} \sqrt{ \gg_{\gamma(t)}(\dot \g(t))}~dt.
\eqnn
The {\it Carnot-Caratheodory distance} (or sub-Riemannian distance) on $M$  associated with 
${\mathcal S}$ is defined as
\bqn\nonumber
d(q_0,q_1)=\inf \{\ell(\g)\mid \g(0)=q_0,\g(T)=q_1, \g\ \mathrm{admissible}\}.
\eqn

{ 
It is a standard fact that $\ell(\g)$ is invariant under reparameterization of the curve $\g$. 
Moreover, if an admissible curve $\g$ minimizes the so-called {\it energy functional} 
$
E(\g)=\int_0^T \gg_{\gamma(t)}(\dot \g(t))~dt
$
with $T$ fixed (and fixed initial and final point)
then $v=\sqrt{\gg_{\gamma(t)}(\dot \g(t))}$ is constant and 
$\g$ is also a minimizer of $\ell(\cdot)$. 
On the other hand a minimizer $\g$ of $\ell(\cdot)$ such that  $v$ is constant is a minimizer of $E(\cdot)$ with $T=\ell(\g)/v$.

The finiteness and the continuity of $d(\cdot,\cdot)$ with respect 
to the topology of $M$ are guaranteed by  the Lie bracket generating 
assumption on the 2-ARS (see \cite{agra-book}).  
The Carnot-Caratheodory distance  endows $M$ with the 
structure of metric space compatible with the topology of $M$ as differential manifold.

%Given a $(n,k)$-\ar\ ${\cal S}$, 

%Notice that the problem of finding a curve minimizing the energy between two fixed points  $q_0,q_1\in M$ is 
%naturally formulated as the optimal control problem
%\bqn
%&&\dot q=\sum_{i=1}^2 u_i X^\mu_i(q)\,,~~~u_i\in\R\,,~~~\mu\in I(q)=\{\mu\in I\mid q\in\Omega^\mu\},\llabel{e-dyn}\\
%&&\int_0^T 
%\sum_{i=1}^2 u_i^2(t)~dt\to \min,~~q(0)=q_0,~~~q(T)=q_1.\llabel{e-cost}
%\eqn
%Here $\mu,u_1,u_2$ are seen as controls and $T$ is fixed. 
%%When $T$ is fixed in such a way that the minimizer is parametrized by arclength,

%It is a standard fact that this optimal control problem  is equivalent to the minimum time problem with 
%controls $u_1,u_2$ satisfying $u_1^2+u_2^2\leq 1$.

When  the 2-ARS is trivializable, the problem of finding a curve minimizing the energy between two fixed points  $q_0,q_1\in M$ is naturally formulated as the distributional optimal control problem  with quadratic cost and fixed final time
\bqn
\dot q=\sum_{i=1}^2 u_i X_i(q)\,,~~~u_i\in\R\,,
~~~\int_0^T 
\sum_{i=1}^2 u_i^2(t)~dt\to\min,~~q(0)=q_0,~~~q(T)=q_1.\eqnn
where $\{X_1,X_2\}$ is an orthonormal frame.
}

\subsection{Geodesics}
A {\it geodesic} for  ${\cal S}$  is a 
curve $\g:[0,T]\to M$ such that 
for every sufficiently small nontrivial interval 
$[t_1,t_2]\subset [0,T]$, $\g|_{[t_1,t_2]}$ is a minimizer of $E(\cdot)$. 
A geodesic for which $\gg_{\gamma(t)}(\dot \g(t))$ is (constantly) 
equal to one is said to be parameterized by arclength. 
The local existence of minimizing geodesics 
is a standard consequence of Filippov Theorem 
(see for instance \cite{agra-book}).  When $M$ is compact any two points of $M$ are connected by a minimizing geodesic.

Locally, in an open set $\Omega$, if  $\{X_1,X_2\}$ is an orthonormal frame, a curve parameterized by arclength is a geodesic if and only if it is the projection on $\Omega$ of a solution of the Hamiltonian system corresponding to the Hamiltonian
 \bqn
\label{eq-HH}
H(q,\lam)=\frac12( (\lam X_1(q))^2 + (\lam X_2(q))^2),~~q\in\Omega,~\lam\in T^\ast_q\Omega.
\eqn
lying on the level set  $H=1/2$.
This is the Pontryagin Maximum Principle \cite{pontryagin-book} in the case of 2-ARSs. Its simple form follows from the absence of abnormal extremals in 2-ARSs, as a consequence of the H\"ormander condition see \cite{ABS}.
 Notice  that when looking for a geodesic $\gamma$  minimizing the  energy from  a submanifold ${\cal T}$ (possibly of dimension zero), one should add the transversality condition $\lam(0)T_{\gamma(0)}{\cal T}=0$.

\subsection{Generic 2-ARSs and normal forms}
\label{s-generic}
A property $(P)$ defined for 2-ARSs  is said to be {\it generic}
if for every rank-2 vector bundle $E$   over $M$, $(P)$ holds for every $\f$ in an  open and dense  subset of the set of  morphisms of vector bundles from $E$ to $TM$, endowed with the 
$C^\infty$-Whitney topology. 

Define  $\bD_1 = \bD$ and $\bD_{k+1}=\bD_k+[\bD,\bD_k]$. We say that $\cal S$ {\it satisfies condition} \HH\  if the following properties hold:
{\bf (i)} $\Zz$ is an
embedded one-dimensional 
submanifold of
$M$;
{\bf (ii)} the points $q\in M$ at which $\bD_2(q)$ is
one-dimensional are isolated;
{\bf (iii)}  $\bD_3(q)=T_qM$ for every $q\in M$. % This is an essential hypothesis needed to prove remarkable results in the context of two-dimensional almost-Riemannian geometry (see
%\cite{euler,ABS} and Theorem \ref{lip-eq} below). Since 
It is not difficult to prove that property \HH\ is generic among 2-ARSs (see  \cite{ABS}). This hypothesis was essential to show Gauss--Bonnet type results for ARSs in \cite{ABS,euler,high-order}. The following  theorem recalls the  local normal forms for ARSs satisfying hypothesis \HH\ (see Figure 1).
%\ref{fig-prenormal}).
\begin{theorem}[\cite{ABS}]
\label{t-normal}
Consider a 2-ARS satisfies \HH. Then
 for every point $q\in M$ there exist a neighborhood $U$ of $q$ and 
 an orthonormal frame $\{X_1,X_2\}$ of the ARS on $U$ such that, up to a change of coordinates, 
 $q=(0,0)$ and $\{X_1,X_2\}$ has one of the forms
$$
\begin{array}{lll}
(\F1) & X_1(x,y)=\frp{}{x}, & X_2(x,y)=e^{\phi(x,y)}\frp{}{y},\\
(\F2) & X_1(x,y)=\frp{}{x}, & X_2(x,y)=xe^{\phi(x,y)}\frp{}{y},\\
(\F3) & X_1(x,y)=\frp{}{x}, & X_2(x,y)=(y-x^2\psi(x))e^{\xi(x,y)}\frp{}{y},
\end{array}
$$
where $\phi$, $\psi$ and $\xi$ are smooth functions such that $\phi(0,y)=0$ and $\psi(0)\neq0$.
\et

Let ${\mathcal S}$ be a 2-ARS satisfying \HH.
 A point $q\in M$ is said to be an
{\it ordinary point} if $\bD(q)=T_q M$, hence, if ${\mathcal
S}$ is locally described by (F1). We call $q$ a  {\it Grushin
point} if $\bD(q)$ is one-dimensional and $\bD_2(q)=T_q M$, i.e., if
the local description (F2) applies. Finally, if
\ugo{$\bD(q)=\bD_2(q)$} has dimension one and $\bD_3(q)=T_q M$
then we say that $q$ is a {\it tangency point} and ${\mathcal
S}$ can be described near $q$ by the normal form (F3). 
%Assume ${\cal S}$ and $M$ to be oriented. Thanks to the hypothesis \HH, $M\setminus \Zz$  splits into two open sets
%$M^+$ and $M^-$ such that $\f:E|_{M^+}\rightarrow TM^+$ is an orientation-preserving isomorphism and $\f:E|_{M^-}\rightarrow TM^-$ is an orientation-reversing isomorphism.
%COPIARE DA LIP EQ. + dare la topologia globale e citare Lip eq.

Notice that under hypotheses {\bf HA} and {\bf HB} of Theorem \ref{t-1}, {\bf H0} is fulfilled and there are no Tangency points. 

In the compact case, the following proposition permits to extend the normal form (F2) to a neighborhood of a connected component of $\Zz$, when there are no tangency points. 
\bp 
\label{p-globalizziagite}
 Consider a 2-ARS on a compact  orientable manifold satisfying ${\bf (H0)}$. Let $\comp$ be a connected component of $\Zz$ containing no tangency points. Then there exists a tubular neighborhood $U$ of $W$  and 
 an orthonormal frame $\{X_1,X_2\}$ of the 2-ARS on $U$ such that, up to a change of coordinates, 
 $W=\{(0,y),~y\in S^1\}$ and $\{X_1,X_2\}$ has the form
 \bqn
\label{global-grushin}
 X_1(x,y)=\frp{}{x}, & X_2(x,y)=x e^{\phi(x,y)}\frp{}{y}.
\eqn
\ep
\proof The proof consists in using as first coordinates the distance from $W$ and 
it is very similar to the proof of Theorem \ref{t-normal}, which contains a local version of Proposition \ref{p-globalizziagite} (see \cite{ABS} pp. 813--814,  proof of Lemma 1 and proof of Theorem 1 for Grushin points). Here we just explain which modifications are necessary.

Since $M$ is assumed to be compact then $W$ is diffeomorphic to $S^1$. Moreover, since $M$ is assumed to be orientable, a  sufficiently small open  tubular neighborhood of $W$ is diffeomorphic to $(-1,1)\times S^1$. 
Then consider a smooth regular parametrization $S^1\ni\al\mapsto w(\al)$ of
$W$.

%{\tt vedere se $-\eps,\eps)$}

Let $\al\mapsto \lam_0(\al)\in T^\ast_{w(\al)}M$ be a smooth map
satisfying  $H(\lam_0(\al),w(\al))=1/2$, where   $H$ is the Hamiltonian of the PMP \r{eq-HH} and
$\lam_0(\al)\perp T_{w(\al)}W$. Notice that such a map exists since we are assuming that $M$ is orientable.

Let $E(t,\al)$ be the solution at time $t$ of the Hamiltonian system given
by the Pontryagin Maximum Principle with initial condition
$(q(0),\lam(0))=(w(\al),\lam_0(\al))$. With the same arguments given in the proof of Lemma 1 in \cite{ABS}, one shows that $E(t,\al)$ is a local diffeomorphism around every point of the type 
$(0,\bar\al)$, $\bar\al\in S^1$. 
Using the fact that $E$ is a global diffeomorphism from $\{0\}\times S^1$ to $W$ and  by suitably reducing $\eps$ one gets that it is a  diffeomorphism between $(-\eps,\eps)\times S^1$  (for some $\eps>0$) and a tubular neighborhood of $W$. 
This permits to use as coordinates in $U$ the pair $(t,\al)$. Indeed $t$ is  the distance from $W$ which we proved to be smooth in $U$. As in \cite{ABS}, one builds an orthonormal frame in $U$ defining the vector field $X_1$ by
\bqn
X_1(t,\al)=\partial_t E(t,\al)=\partial_t.\nonumber
\eqn
As in \cite{ABS} one obtains that the second vector field of the orthonormal frame has to be of the form  $t e^{\phi(t,\al)}\frp{}{\al}$. $\square$

\brem\label{r-or}
Notice that in Proposition \ref{p-globalizziagite}, the hypothesis that $M$ is orientable can be weakened by requiring that  $W$ admits an open  tubular neighborhood diffeomorphic to $(-1,1)\times S^1$
\erem
 
%\brem
% Almost-Riemannian structures attracted a certain interest in the last years. In \cite{ABS,high-order} a Gauss--Bonnet-type formula was obtained in the case without tangency points and under generic conditions. This formula was generalized in \cite{euler} in the more intricate case in which tangency points are present.   In \cite{BCGS} a necessary and sufficient condition for two  2-ARSs  on the same compact manifold $M$ to be   Lipschitz equivalent was given. This equivalence was established  in terms of graphs associated with the structures. 
% 
%Tangency points are the most difficult to handle due to the fact that the asymptotic of the distance is different from the two sides of the singular set.  
%In \cite{BCGJ} the authors gave a description of the geometry of the nilpotent approximation at a tangency point, provided jets of the exponential map and a description of the cut and conjugate loci from a tangency point in the  generic case. In \cite{kresta}, the authors considered the problem of finding a completely reduced normal form at a generic tangency point.
%\erem 

\section{The Grushin case}
\label{s-grushin}
As already mentioned the Grushin plane  
is the trivializable almost Riemannian metric on the $(x,y)$ plane 
for which an orthonormal basis is given by 
\bqn
X_1=
\left(\ba{c}
1\\0
\ea
\right),~~~X_2=
\left(\ba{c}
0\\x
\ea
\right).
\eqn
In the sense of Definition \ref{d-ARS} it can be seen as a triple $(E, \f,\langle\cdot,\cdot\rangle)$ where $E=\R^2\times\R^2$, $\f((x,y),(a,b))=((x,y),(a, b x))$ and $\langle\cdot,\cdot\rangle$ is the standard Euclidean \ugo{metric}.

\subsection{Geodesics of the Grushin plane}
\label{ss-geo-grushin}

 In this section we recall how to compute the geodesics for the Grushin plane, with the purpose of stressing  that they can cross the singular set with no singularities. 

Setting $q=(x,y)$ and $\lambda=(\lambda_1,\lambda_2)$, the Hamiltonian \r{eq-HH} is given by
\bqn
\label{h-g}
H(x,y,\lambda_1,\lambda_2)=\frac12(\lambda_1^2+x^2\lambda_2^2)
\eqn
and the corresponding Hamiltonian equations are:
\bqn
&&\dot x=\lambda_1,~~~~\dot \lambda_1=-x \lambda_2^2\nn\\
&&\dot y=x^2\lambda_2,~~\dot \lambda_2=0 
\eqn

 Geodesics parameterized by arclength  are projections on the $(x,y)$ plane of solutions of these equations, lying on the level set  $H=1/2$.
We study the geodesics starting from {\bf i)} a Grushin point, e.g. $(0,0)$ {\bf ii)} an ordinary point, e.g. $(-1,0)$.

\medskip
\noindent
\underline{{\bf Case} $(x(0),y(0))=(0,0)$} \\
In this case the condition $H(x(0),y(0),\lambda_1(0),\lambda_2(0))=1/2$ implies that we have two families of geodesics corresponding respectively to
$\lambda_1(0)=\pm1,~~~\lambda_2(0)=:a\in\R$. 
Their expression can be easily obtained and it is given by:

\bqn
\left\{
\ba{lll}
x(t)=\pm t,& y(t)=0 &\mbox{ if } a=0\\
x(t)=\pm \frac{\sin(a t)}{a},&y(t)=\frac{2 a t - \sin(2 a t)}{4 a^2}&\mbox{ if } a\neq0
\ea\right.
\eqn
{   Some geodesics are plotted in Figure 2 
%\ref{f-grushin1} 
together with the ``front'' i.e. the end point of all geodesics at time $t=1$. Notice that geodesics start horizontally. The particular form of the front shows the presence of a conjugate locus accumulating to the origin.}

\begin{figure}
\label{f-grushin1}
\begin{center}
~\includegraphics[width=8truecm]{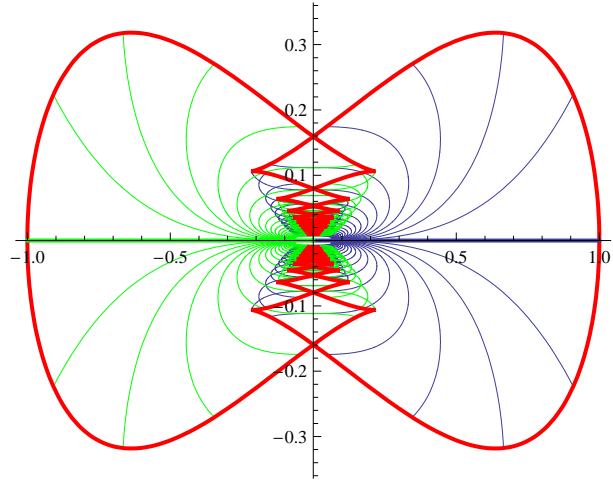}\\
\caption{Geodesics and front for the Grushin plane, starting from the singular set.}
\end{center}
\end{figure}

\medskip
\noindent
\underline{{\bf Case} $(x(0),y(0))=(-1,0)$}\\
 In this case the condition $H(x(0),y(0),\lambda_1(0),\lambda_2(0))=1/2$ becomes $\lambda_1^2+\lambda_2^2=1$
 and it is convenient to set $\lambda_1=\cos(\th), \lambda_2=\sin(\th),~~ \th\in S^1.$
The expression of the geodesics is given by:
\bqn
\left\{
\ba{lll}
x(t)= t-1,& y(t)=0 &\mbox{ if } \th=0\\
x(t)=-t-1,& y(t)=0 &\mbox{ if } \th=\pi\\
\displaystyle x(t)=-\frac{\sin (\th-t \sin (\th))}{\sin(\th)} &\displaystyle
y(t)=\frac{2 t - 2 \cos(\th) +\frac{\sin(2 \th - 2 t \sin(\th))}{ \sin(\th) }}{4\sin(\th)}
&\mbox{ if } \th\notin\{ 0,\pi\}
\ea\right.
\eqn
Some geodesics are plotted in Figure 3
 %\ref{f-grushin4.8} 
 together with the ``front''  at time $t=4.8$. Notice that geodesics pass horizontally through $\Zz$, with no singularities. The particular form of the front shows the presence of a conjugate locus. Geodesics can have conjugate times only after intersecting $\Zz$. Before it is impossible since they are Riemannian and the curvature is negative. 

\begin{figure}[h!]
\label{f-grushin4.8}
\begin{center}
~\includegraphics[width=8truecm]{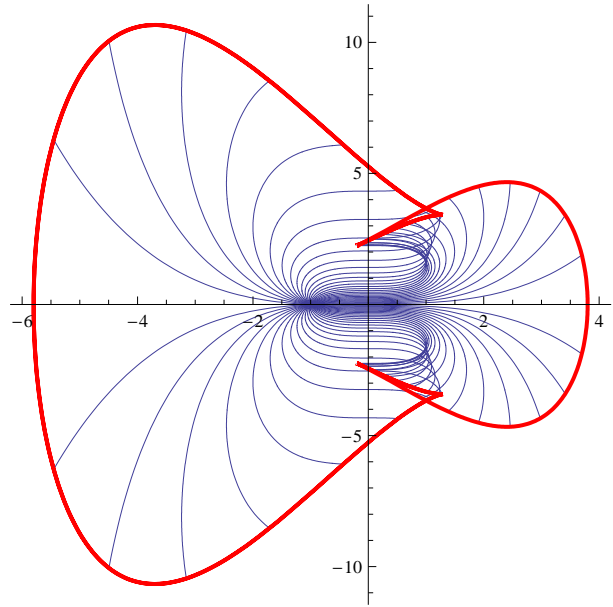}
\caption{Geodesics and front for the Grushin plane, starting from a Riemannian point}
\end{center}
\end{figure}

 \subsection{The Laplace-Beltrami operator on the Grushin plane}
In this section we illustrate in details the results of the paper in the Grushin case.  For the study of the Laplace Beltrami operator it will be convenient to compactify the $y$ direction considering $x\in\R$ and $y\in \Tu$.

As mentioned in the introduction, the Laplace-Beltrami operator on the Grushin plane can be written as 
$$
\Delta=\partial_{x}^2+x^2\partial_{y}^2-\frac{1}{x}\partial_{x}.
$$
We expect that the Laplacian is the Friedrichs extension associated to the positive quadratic form
\bna
Q(u,u)&=&\int_{\R^2}|\nabla u|^2~d\omega =\int_{\R^2}\frac{1}{|x|} |\partial_{x} u|^2+|x|\,|\partial_{y} u|^2dxdy
\ena
Let us make the change of variable $f=\sqrt{|x|}g$ which is unitary from $L^2(\R\times \Tu,d\omega)$ to $L^2(\R\times \Tu,dx\,dy)$, so that $(f_1,f_2)_{L^2(\R\times \Tu,d\omega)}= (g_1,g_2)_{L^2(\R\times \Tu,dx\,dy)}$.

We compute the operator in the new variable:
\bna
\Delta f&=&\partial_x^2f+x^2\partial_y^2f-\frac{1}{|x|}\partial_xf\\
&=&
\sqrt{|x|}\left[ \partial_x^2 g+x^2 \partial_y^2g -\frac{3}{4x^2}g\right]=:\sqrt{|x|}Lg.
\ena
%Et on a bien alors (modulo le produit scalaire sur la deuxième composante)
%\bna
%(\Delta_{Gm}g_1,g_2){L^2_{\R^2}}=-(\nabla_{Gm}g_1,\nabla_{Gm}g_2){L^2_{\R^2}}=(\Delta_{Gm}g_1,g_2){L^2_{\R^2}}
%\ena
Hence we are left to study the operator $Lg=\partial_x^2 g+x^2 \partial_y^2g -\frac{3}{4x^2}g$ on $L^2(\R\times \Tu,dx\,dy)$. Decomposing a function $g(x,y)$ in the Fourier basis in the variable $y$, we get the decomposition 
\bna
L^2(\R\times \Tu)\approx \oplus_{k\in \Z}^{\perp} H_k
\ena
where $H_k=L^2(\R)$ and the operator $L$ acts on each $H_k$ by 
\bna
(L u)_k(x)=\partial_x^2 u_k-x^2k^2 u_k -\frac{3}{4x^2} u_k=\partial_x^2 u_k-V_k u_k.
\ena
with $V_k(x)\geq \frac{3}{4x^2}$. But, we know that in dimension $1$, the operator $-\partial_x^2 +V$ with domain $C^{\infty}_0(]0,+\infty[)$ is essentially self-adjoint on $L^2(]0,+\infty[)$ if $V\geq \frac{3}{4x^2}$ (see \cite{ReedSimon}, Theorem X.10 for the proof of the limit point case at $0$ and Theorem X.8 at $+\infty$). Hence, each operator $\partial_x^2 u_k-V_k$ is essentially self-adjoint. As a consequence $L$ is essentially self-adjoint as well.

\brem
\label{r-alfa}
 $V=\frac{3}{4x^2}$ is exactly the limit singular potential to have  $-\partial_x^2+V$  in the limit point case at zero and hence essentially self-adjoint (see \cite[Theorem X.10, page 159]{ReedSimon}). It is  surprising to find the same constant $\frac{3}{4}$ for the Grushin operator.

Notice that if we consider as orthonormal basis of the metric the vector fields  $X=\begin{pmatrix}1\\0\end{pmatrix}$ and $Y=\begin{pmatrix}0\\|x|^{\alpha}\end{pmatrix}$ with $\alpha>0$, we obtain  for the Laplace-Beltrami operator,
\bna
\Delta u=\partial_x^2u+|x|^{2\alpha}\partial_y^2u-\frac{\alpha}{x}\partial_xu.
\ena
Making the change of variable $u=|x|^{\alpha/2}v$ (which is unitary from $L^2(\R\times\Tu,d\omega)$ to $L^2(\R\times\Tu,dx\,dy)$) we get for the transformed operator
\bna
\widetilde{\Delta} v=\partial_x^2v-\frac{\alpha}{2}\left(\frac{\alpha}{2}+1 \right) \frac{1}{x^2}+|x|^{2\alpha}\partial_y^2v,
\ena
which is essentially self-adjoint if and only if $\alpha \in
 %]-\infty,-3]\cup 
 [1,+\infty[$.
\erem

\subsection{$\sqrt[1/\al]{\mbox{Grushin}}$}
\label{s-sqrt}
In this section, motivated by Remark \ref{r-alfa}, we briefly discuss the generalized Riemannian structures on $\R^2$ for which an orthonormal basis is given by 
\bqn
X_1=\left(\ba{c}1\\0 \ea\right),~~~X_2=\left(\ba{c}0\\|x|^\al \ea\right),~~\al\in(0,1).
\eqn
Notice that this generalized-Riemannian structure does not fit the definition of almost-Riemannian structures, since the vector fields are not smooth. As usual call $\Zz$ the set in which the two vector fields are not linearly independent. It is interesting to notice the following facts, which in some cases are very different with respect to the standard Grushin case:
\bi
\iii The Lie bracket between $X_1$ and $X_2$ is not well defined on $\Zz$. Indeed $[X_1,X_2]=(0,\al |x|^{(\al-1)} )^T$. However it is not difficult to prove that the corresponding distance is continuous and endows  $\R^2$ with its standard topology (note that it is not the case if $\alpha>0$).
\iii The curvature is given by $K=-\frac{\al(1+\al)}{x^2}$ and the area element  by $d\omega=\frac{1}{|x|^\al}dx\, dy$. Hence the area element explodes, but the area of a 
bounded open set intersecting $\Zz$  is finite (and its diameter is finite as well). 
This is certainly the reason why the Laplace-Beltrami operator is not essentially self-adjoint and, as a consequence, a quantum particle or the heat flow can pass through the set $\Zz$.
Notice that for the Grushin plane a bounded open set intersecting $\Zz$ has finite diameter but infinite area.

\iii There exists regular curves passing through $\Zz$ with a velocity not belonging to the span of $X_1$ and $X_2$ and having finite length. For instance take $\al=1/2$ and the curve $\gamma(t):=(t,t)$, $t\in[-1,1]$. We have  $\dot\gamma(0)\notin$~span$\{X_1(\gamma(0)),X_2(\gamma(0)\}$ but its length is
finite. Indeed it corresponds to the controls $u_1=1$ and $u_2=1/|t|^{1/2}$. Hence $\ell(\gamma)=\int_{-1}^{1}\sqrt{1+|t|^{-1}}dt=2 (\sqrt{2} + \text{arcsinh}(1))\sim4.59.$ Such phenomenon does not appear on the singular set of a  2-ARS, as a consequence of the smoothness of  orthonormal frames.

\ei

\section{Laplace-Beltrami on an Almost-Riemannian structure}
\label{s-camillo}
\subsection{Computation in local coordinates}
The aim of this subsection is the proof of the following theorem.
\begin{theorem}
\label{Grushinvaress}
 Let $g$ a Riemannian metric on $\Omega=\left\{x=(x,y)\in \R \times \Tu  \left| x\neq 0\right.\right\}$ of the form diag$\left(1, x^{-2} e^{-2\Phi(x,y)}\right)$ where  $\Phi$ is a smooth function which is constant for large $|x|$. Let $\Delta$ be the corresponding Laplace-Beltrami operator with domain $C^{\infty}_0(\Omega)$.

Then,  $\Delta$ is essentially self-adjoint for $L^2(\R\times \Tu,d\omega)$, where $d\omega$ is the corresponding  Riemannian area on $\Omega$.  

The same result holds by replacing $\Omega$ with $\Omega_{\pm}$ where  $\Omega_-=\left\{x=(x,y)\in \R \times \Tu  \left| x< 0\right.\right\}$  and $\Omega_+=\left\{x=(x,y)\in \R \times \Tu  \left| x> 0\right.\right\}$.

\end{theorem}
A simple consequence of this theorem is that the only self-adjoint Laplace operator that can be constructed by extension of the natural one (i.e. the Laplace-Beltrami operator defined on $C^{\infty}_0(\Omega)$) preserves the decomposition  $L^2(\Omega,d\omega)=L^2(\Omega_-,d\omega)\oplus ^{\perp} L^2(\Omega_+,d\omega)$ . 

Indeed, if we denote $\Delta^-$ (resp $\Delta^+$) the unique self-adjoint extension of the natural Laplace operator with domain $C^{\infty}_0(\Omega_-)$ (resp $C^{\infty}_0(\Omega_+)$), we can define  the operator $ \bar{\Delta} (u_-+u_+ )=\Delta^- u_-+ \Delta^+u_+$ with the natural domain inherited from the one of $\Delta^-$ and $ \Delta^+$.   The operator $\bar{\Delta}$ is a self-adjoint extension of the Laplace-Beltrami operator with domain $C^{\infty}_0(\Omega)$ and it is the unique one by Theorem \ref{Grushinvaress}.  Another way of seeing this fact is from the point of view of the evolution equation. 
\bc
With the notations of Theorem \ref{Grushinvaress}, consider the unique solution $u$ of the Schrödinger equation (according to the self-adjoint extension defined in the previous theorem),
\bneqn
i\partial_t u+\Delta u&=&0\\
u(0)&=&u_0\in  L^2(\R\times \Tu,d\omega)
\eneqn
with $u_0$ supported in $\Omega_+$. Then, $u(t)$ is supported in $\Omega_+$  for any $t\geq 0$.
The same holds for the solution of the heat or for the solution of the wave equation. 
\ec
This corollary is a simple consequence of the essential self-adjointness of the Laplace-Beltrami operator as discussed above.

This is in strong contrast with the classical dynamics associated with this metric where the geodesics can "cross" the  barrier $\left\{x=0\right\}$ (see Figure 3).
%\ref{f-grushin4.8}).
The proof of Theorem \ref{Grushinvaress} relies on a change of variable  similar to the one we  did for the Grushin operator which leads to an operator that can be written as ``Grushin type + singular potential''. Then, we follow the proof of Kalf and Walter \cite{KalfWalter_essential_singular}  (which was himself inspired by B. Simon \cite{Simon_essential_singular}) using Kato inequality for the main part. The other part can be treated by perturbation theory. Note that some related results for singular potential were proven in \cite{DonnellySchrod} and \cite{MaedaSchrodess_affin}.
\bnp[Proof of Theorem \ref{Grushinvaress}]
In coordinates, the Laplace-Beltrami operator  has the form
\bna
\Delta=\frac{1}{\omega}\sum_{i,j}\partial_i\left[\omega g^{ij}\partial_j\right]
\ena
where we have denoted $\omega=\sqrt{\text{det} g}$ the Riemannian volume.

For $u\in L^2(\R\times \Tu,d\omega)$, let us make the change of variable $ u=\omega^{-1/2}v$ which is unitary from $L^2(\R\times \Tu,d\omega)$ to $L^2(\R\times \Tu,dx\,dy)$ and let us  compute its action in the new variable.
\bna
&&\omega^{1/2}\Delta\omega^{-1/2}v:=\widetilde{\Delta}v\\&=&\frac{1}{\omega^{1/2}}\sum_{i,j}\partial_i\left[\omega g^{ij}\partial_j(\omega^{-1/2}v)\right]\\
&=&\frac{1}{\omega^{1/2}}\sum_{i,j}\partial_i\left[\omega g^{ij}(\partial_j\omega^{-1/2})v+g^{ij}(\omega^{1/2})\partial_jv\right]\\
&=&\frac{1}{\omega^{1/2}}\sum_{i,j}\left[ \partial_i\left[g^{ij}\omega(\partial_j\omega^{-1/2})\right]v+ g^{ij}\omega(\partial_j\omega^{-1/2})\partial_iv+(\partial_ig^{ij})\omega^{1/2}\partial_jv+g^{ij}(\partial_i\omega^{1/2})\partial_jv+g^{ij}(\omega^{1/2})\partial_{ij}v\right]\\
&=&\frac{1}{\omega^{1/2}}\sum_{i,j}\left[ \partial_i\left[g^{ij}\omega(\partial_j\omega^{-1/2})\right]v+ g^{ij}\omega(\partial_i\omega^{-1/2})\partial_jv+(\partial_ig^{ij})\omega^{1/2}\partial_jv+g^{ij}(\partial_i\omega^{1/2})\partial_jv+g^{ij}(\omega^{1/2})\partial_{ij}v\right]
\ena
where we have used the symmetry of $g$.

We have:  $\omega(\partial_i\omega^{-1/2})+(\partial_i\omega^{1/2})= -\frac{1}{2}\omega \omega^{-3/2}\partial_i\omega+\frac{1}{2} \omega^{-1/2} \partial_i\omega=0$. 
Moreover let 
\bqn
\mu:=\sum_{i,j}\frac{1}{\omega^{1/2}}\partial_i\left[g^{ij}\omega(\partial_j\omega^{-1/2})\right]
&=&-\sum_{i,j}\frac{1}{2\omega^{1/2}}\partial_i\left[g^{ij}(\partial_j\omega)\omega^{-1/2})\right]\nn\\
&=&\sum_{i,j}\left[- \frac{1}{2}\partial_i(g^{ij})\partial_j(\ln\omega)-\frac{1}{2\omega}g^{ij}\partial_{ij}\omega+\frac{1}{4\omega^2}g^{ij}(\partial_i \omega)(\partial_j \omega)\right].\nn
\eqn
Hence we get 
\bna
\widetilde{\Delta}v&=& \sum_{i,j}g^{ij}\partial_{ij}v+(\partial_ig^{ij})\partial_jv+\mu v\\
&=&div_{\eucl} (\nabla v)+\mu v\nn\\
&=& div_{\eucl} (G^{-1}\nabla_{\eucl} v)+\mu v.
\ena
In the last formula we have used matrix notation where $G^{-1}=(g^{ij})$ and $\div_{\eucl}$, $\nabla_{\eucl}$ are the divergence and the gradient of the Euclidian space. Now let us specify the computation in our diagonal metric.
\bna
\omega&=&\frac{1}{|x|}e^{-\Phi},\\
\partial_1\omega&=&-\frac{\sgn(x)}{|x|^2}e^{-\Phi}-\frac{\partial_1 \Phi}{|x|}e^{-\Phi},\\
\partial_2\omega&=&-\frac{\partial_2 \Phi}{|x|}e^{-\Phi}; \quad \partial_2^2\omega=\frac{(\partial_2 \Phi)^2}{|x|}e^{-\Phi}-\frac{\partial_2^2 \Phi}{|x|}e^{-\Phi},\\
\partial_1^2\omega&=&\frac{2}{|x|^3}e^{-\Phi}+\frac{2\sgn(x)\partial_1 \Phi}{|x|^2}e^{-\Phi}+\frac{(\partial_1 \Phi)^2}{|x|}e^{-\Phi}-\frac{\partial_1^2 \Phi}{|x|}e^{-\Phi},\\
\mu&=&-\frac{1}{|x|^2}-\frac{\sgn(x)\partial_1 \Phi}{|x|}-\frac{(\partial_1 \Phi)^2}{2}+\frac{\partial_1^2 \Phi}{2}\\
&&+\frac{1}{4|x|^2}+\frac{\sgn(x)\partial_1 \Phi}{2|x|}+\frac{(\partial_1 \Phi)^2}{4}\\
&&+x^2e^{2\Phi}\left[ (\partial_2\Phi)^2+\frac{1}{2}\left[\partial_2^2 \Phi-(\partial_2 \Phi)^2\right]+\frac{1}{4}(\partial_2 \Phi)^2\right]\\
&=&-\frac{3}{4|x|^2}-\frac{\sgn(x)\partial_1 \Phi}{2|x|}-\frac{(\partial_1 \Phi)^2}{4}+\frac{\partial_1^2 \Phi}{2}\\
&&+x^2e^{2\Phi}\left[\frac{1}{2}\partial_2^2 \Phi+\frac{3}{4}(\partial_2 \Phi)^2\right]\\
&:=&-\frac{3}{4|x|^2}-\mu_2.%-\frac{1}{2}(n-1)(n-2)|x|^{-2}+\grando{|x|^{n-2}+|x|^{n-1}}\\
%&&+\frac{1}{4}(n-1)^2|x|^{-2}+\grando{|x|^{n-2}+|x|^{n-1}}
\ena

Now, we are able to conclude by the following two Lemmas.  The first one proves that the main part without $\mu_2$ is essentially self-adjoint while the second one treats $\mu_2$ as a perturbation.
\bl
\label{lemmeessauto}
Let $G^{-1}$  given by diag$\left(1, x^2 e^{2\Phi(x,y)}\right)$ with $|\Phi| \leq C$ for $C>0$.

The operator $L: v \mapsto -\div_{\eucl} (G^{-1}\nabla_{\eucl} v) +\frac{3}{4|x|^2}v$ defined on $L^2(\R\times \Tu,dx\,dy)$ with domain $C^{\infty}_0(\Omega)$ is positive and essentially self-adjoint. 
\el
Hence, since $L$ defines a positive symmetric operator, its unique self-adjoint extension is the Friedrichs extension
(in the following still denoted by $L$).
\bl
\label{lemmeperturb}
The operator of multiplication by $\mu_2$ (defined with domain $C^{\infty}_0(\Omega))$) is infinitesimally small with respect to $L$. So, by the theorem of Kato-Rellich (Theorem X.12 of \cite{ReedSimon}), $L+\mu_2$ is essentially self-adjoint.
\el
To conclude the proof of Theorem \ref{Grushinvaress} we are left to prove the two lemmas.

\bnp[Proof of Lemma \ref{lemmeessauto}]
We follow closely the proof of Kalf and Walter \cite{KalfWalter_essential_singular}.  Actually, we mimick the proof of the 1-D case and we use the fact that our operator is the classical Laplacian with potential for function only depending on $x$: $-\partial_1^2+\frac{3}{4|x|^2}$.
It is enough to prove that $\textnormal{Ran}(L+I)$ is dense in $L^2$, see Theorem X.26 of \cite{ReedSimon}. So, let $h\in L^2(\Omega)$ such that $h \perp \textnormal{Ran}(L+I)$ and let us prove $h=0$ (we can assume $h$ real valued without loss of generality). 

Set $\lettre(x)=\frac{|x|^{3/2}e^{-|x|}}{1+|x|^{3/2}}$. We have $\lettre(x)\geq 0$ and $\lettre \in L^2(\Omega)$. We obtain, for $x\neq 0$,
\bna
-\div_\eucl (G^{-1}\nabla_\eucl \lettre) +\frac{3}{4|x|^2} \lettre+\lettre&=&-\partial_1^2\lettre+\frac{3}{4|x|^2}\lettre+\lettre\\
&=&\frac{3|x|^{1/2}e^{-|x|}}{4(1+|x|^{3/2})^3}(|x|^2+4|x|^{3/2}+7|x|^{1/2}+4):=X(x)
\ena 
%(j'ai fait le calcul que pour $x>0$)
Let $\eta$, $\zeta \in C^{\infty}([0,+\infty[)$ so that 
\bna
\eta(s)&=&0 \textnormal{ if } s\leq 1/2\\
&=&1\textnormal{ if } s\geq 1\\
\zeta(s)&=&1 \textnormal{ if } s\leq 1\\
&=&0\textnormal{ if } s\geq 2
\ena
Denote $\lettre_n(x) =\eta(n |x|)\zeta(\frac{|x|}{n})\lettre(x)\in C^{\infty}_0(\Omega)$ so that $L\lettre_n$  is well defined in $L^2(\Omega)$. Moreover, we have
\bna
\lettre(x)=\grando{|x|^{3/2}}, &\quad &|\nabla \lettre(x)|=\grando{|x|^{1/2}}, \textnormal{ as }x\rightarrow 0\\
\lettre(x)=\grando{e^{-|x|}}, &\quad & |\nabla \lettre(x)|=\grando{e^{-|x|}}, \textnormal{ as }x\rightarrow \infty
\ena
 so, $L\lettre_n+\lettre_n=\eta(n |x|)\zeta(\frac{|x|}{n}) X-2 \partial_1\left[\eta(n |x|)\zeta(\frac{|x|}{n}) \right]\partial_1\lettre-\partial_1^2\left[\eta(n |x|)\zeta(\frac{|x|}{n})\right]\lettre $ converges weakly in $L^2$ 
 %\bleu{(quelques calculs de commutateurs à vérifier)} 
 to $X\in L^2(\Omega)$ such that $X(x)\geq 0$. So, we can write
 \bna
 \left(X,|h|\right)_{L^2}=\limvar{n}{\infty} \left(L\lettre_n+\lettre_n,|h|\right)_{L^2}\geq 0
 \ena

But, we will use the Kato inequality: 
\bl[Kato's inequality, see Lemma A of \cite{Kato_ineg}]
\label{lemmeKato}
~\\
For $v\in L^1_{loc}(\Omega)$ real valued such that $div_{\eucl} (G^{-1}\nabla_{\eucl} v)\in L^1_{loc}(\Omega)$, we have, in the sense of distributions on $\Omega$
\bna
-\div_{\eucl} (G^{-1}\nabla_{\eucl} |v|)\leq -(\sgn v)\div_{\eucl} (G^{-1}\nabla_{\eucl} v).
\ena
\el
Notice that the Kato inequality of \cite{Kato_ineg} can only be applied on $\Omega_\eps:=\{(x,y)\in\Omega ~:~|x|>\eps\}$ where the metric is positive definite. This does not create difficulties since we are applying this Lemma in the sense of distributions on $\Omega$.

Applying the Kato Lemma to $h$ and $\lettre_n>0$ as test function, we get
\bna
\left(L\lettre_n+\lettre_n,|h|\right)_{L^2}\leq \left(\lettre_n,(\sgn h)(Lh+h)\right)_{L^2}\leq 0.
\ena
So $\left(X,|h|\right)_{L^2}=0$ and $h\equiv 0$.
\enp
\bnp[Proof of Lemma \ref{lemmeperturb}]
The quadratic form associated to $L+1$ is
\bna
Q(u,u)=\int_{\Omega} (G^{-1} \nabla_e u\cdot  \nabla_e u)+\frac{3}{4|x|^2}\left| u\right|^2+\left| u\right|^2~dx\,dy
\ena

For $u\in C^{\infty}_0(\Omega)$, and by Cauchy-Schwarz inequality, we have 
\bna
Q(u,u)= \left( Lu,u\right)_{L^2} +\norL{u}^2\leq \norL{Lu}\norL{u}+\norL{u}^2\leq \varepsilon \norL{Lu}^2+\left( \frac{1}{4\varepsilon}+1\right) \norL{u}^2
\ena

Moreover, since $\Phi$ is smooth and constant for large $x$, we have
\bna
\norL{\mu_2 u}^2 &\leq& C\left( \norL{\frac{u}{|x|}}+\norL{u}^2\right) \\
&\leq&C Q(u,u)\leq C\varepsilon \nor{Lu}{L^2}^2+(C_{\varepsilon}+1)\nor{u}{L^2}^2.
\ena
This yields the result since $\varepsilon$ is arbitrary small and the estimate can be extended to any $u\in D(L)$ by taking closure.
\enp
\noindent
The proof of Theorem \ref{Grushinvaress} is concluded. 
\enp
  {
The next Lemma is the first step to establish the compactness of the resolvent of $\Delta$ in the compact case. Of course, since here we are dealing with the case of $\R\times\Tu$, we can not expect compactness without adding a cut-off function.
\bl
\label{lmcompactcoord}
Denote by $\Delta$ the positive self-adjoint operator defined by Theorem \ref{Grushinvaress}. Then, the truncated resolvent $\rho(x)(-\Delta +1)^{-1}$, where $\rho \in C^{\infty}_0(\R)$, is well-defined and compact from $L^2(\R\times\Tu,d\omega)$ to itself. 
\el}
\bnp
The fact that $(-\Delta +1)^{-1}$ is well defined comes from the positivity of $-\Delta$. For the compactness it is equivalent to prove that the operator $\rho(x)(-\widetilde{\Delta}+1)^{-1}$ is compact when defined on $L^2(\R\times\Tu,dx\,dy)$ with the metric induced by the Lebesgue measure.

We begin by showing that it is the case for $L$ defined above. Let $u_n$ be a real valued sequence in $L^2(\R\times\Tu,dx\,dy)$ with $L^2$ norm bounded by $1$. Denote $\widetilde{f_n}:=\rho(x)(L +1)^{-1}u_n:=\rho(x)f_n $. In particular, $Q(f_n,f_n)=((L+1)f_n,f_n)=(u_n,f_n)=\int_{\Omega} G^{-1}\nabla_{\eucl} f_n\cdot \nabla_{\eucl} f_n+\frac{3}{4|x|^2}\left| f_n\right|^2+\left| f_n\right|^2~dx\,dy$ is bounded in $L^2(\R\times\Tu,dx\,dy)$. Set $\varepsilon>0$. Take $\varphi\in C^{\infty}_0(\R)$, supported in $]-\varepsilon,\varepsilon[$ with $0\leq \varphi \leq 1$ and $\varphi \equiv 1$ in a neighborhood of $0$. Split $\widetilde{f_n}=\varphi(x)\widetilde{f_n}+(1-\varphi(x))\widetilde{f_n}:=f_{n,1}+f_{n,2}$. Then, we have $|x\varphi(x)|^2 \leq\varepsilon^2$ and $\nor{ \varphi f_n}{L^2}^2= \int \frac{|x\varphi(x)|^2 f_n}{|x|^2}\leq \varepsilon^2 \frac{4}{3}Q(f_n,f_n)\leq C\varepsilon^2$. This gives also $\nor{f_{n,1}}{L^2}\leq C\varepsilon$. Now that $\varphi$ is fixed, $f_{n,2}=(1-\varphi(x)) \widetilde{f_n}=(1-\varphi(x))\rho(x)f_n$ is supported in a fixed compact of $\Omega$ and bounded in $H^1$. Indeed, 
\bna
\norL{\nabla_{\eucl} f_{n,2}}^2 \leq 2 (\norL{(1-\varphi(x))\rho(x)\nabla_{\eucl} f_n}^2 +  \norL{\partial_{x}(\varphi(x)\rho(x)) f_n}^2)\leq CQ(f_n,f_n)\leq C
\ena
where the second inequality comes from the fact that there exist a constant $C$ such that $Id \leq C G^{-1}$ (in the sense of quadratic forms) on the support of $(1-\varphi(x))\rho(x)$.
We conclude by invoking the compact embedding of $H^1$ into $L^2$ on compact sets so that $f_{n,2}$ is convergent up to extraction. Actually, we have proved that for any $\varepsilon>0$, we can find an extraction $\gamma$ such that $f_{\gamma(n)}$ can be written $f_{\gamma(n)}=f_{\gamma(n),1}+f_{\gamma(n),2}$ with $\norL{f_{\gamma(n),1}}\leq \varepsilon$ and $f_{\gamma(n),2}$ convergent. By choosing $\varepsilon=\frac{1}{p}$, $p\in \N$, and by a diagonal extraction argument, we easily get that we can find a subsequence such that $f_{\gamma(n)}$ is a Cauchy sequence and converges. This gives that $\rho(L+1)^{-1}$ is compact. It remains to prove the same result for $L+\mu_2$.

Again, let $u_n$ be a sequence in $L^2(\R\times\Tu,dx\,dy)$ with $L^2$ norm bounded by $1$. Denote $\widetilde{f_n}:=\rho(x)(L+\mu_2+1)^{-1}u_n:=\rho(x)f_n $. Thanks to Lemma \ref{lemmeperturb}, we get
\bna
\norL{(L+1)f_n} &\leq &\norL{(L+\mu_2+1)f_n}+\norL{\mu_2 f_n}\leq \norL{u_n}+ \varepsilon \norL{(L+1)f_n}+C_{\varepsilon}\norL{f_n}\\
&\leq &C_{\varepsilon}+ \varepsilon \norL{(L+1)f_n}
\ena
so that we get by absorption for $\varepsilon$ small enough that $(L+1)f_n$ is bounded. Since we have proved that $\rho (L+1)^{-1}$ is a compact operator, we get that $\widetilde{f_n}=\rho(x)f_n=\rho (L+1)^{-1}(L+1)f_n$ is relativelly compact.
\enp
\subsection{Case of a compact manifold}
In this section, using the previous results, we prove the  three statements of Theorem \ref{t-1}. 
 For the first statement, we have to prove the following.
\begin{proposition}
\label{p-thmesscompact}
Consider a 2-ARS on a compact manifold $M$  satisfying hypotheses {\bf HA} and {\bf HB}. Denote by $\Delta$ the Laplace operator defined for functions of $C^{\infty}_0(\Omega)$ where $\Omega=M\setminus \Zz$.  Then, $\Delta$ is essentially self-adjoint on $L^2(M,d\omega)$.
\end{proposition}
\bnp
Let $M=\bigcup_{i\in I} \Omega_i$ be a finite covering of $M$ such that for every connected component $W$ of  $\Zz$ there exists $i\in I$ such that $W\subset \Omega_i$ and $W\cap \Omega_j=\emptyset$ for $j\in I$, $j\neq i$ and an orthonormal frame for the 2-ARS in $\Omega_i$  is given by \r{global-grushin}. Moreover assume that if $\Omega_j$ ($j\in I$) does not contain any Grushin point, then  an orthonormal frame for the 2-ARS in $\Omega_i$  is given by the normal form (F1). This is possible thanks to Theorem \ref{t-normal} and Proposition \ref{p-globalizziagite}.

Theorem \ref{Grushinvaress} yields the result in local coordinates around a connected component of the singular set. We only have to extend it to the whole manifold. 

So, let $(\Psi_i)_{i\in I}$ be a partition of unity associated to $(\Omega_i)_{i\in I}$, that is $\Psi_i\in C^{\infty}_0(\Omega_i)$ and $\sum_{i\in I} \Psi_i =1$. We can also assume that $\Psi_i\equiv 0$ or $\Psi_i \equiv 1$ in a neighborhood of $\Zz$, so that $\nabla \Psi_i$ and $\Delta \Psi_i$ are $C^{\infty}$ compactly supported in $\Omega$. Here we use the fact that we have global coordinates around every connected component of $\Zz$.

So again, let $u\in L^2(\Omega)$ such that $u \perp \textnormal{Ran}(-\Delta+I)$ and let us prove $u=0$ (we can assume $u$ real valued without loss of generality). 
We denote $u_i=\Psi_i u \in L^2(\Omega_i)$.

Since the domain of definition of $\Delta$ contains $C^{\infty}_0(\Omega)$, we get that $u$ is solution of $-\Delta u+u=0$ in the sense of distributions on $\Omega$. By elliptic regularity, we get that $u\in C^{\infty}(\Omega)$. So, the only problem is the possibility to compute integration by part around the degeneracy points.

Let $\Omega_i$ with $i\in I$ be an open set around which we can find a coordinate system so that the metric $g$ takes the diagonal form $\left(1, x^{-2} e^{-2\Phi(x,y)}\right)$ that we  extend  arbitrarily on  $\R \times\Tu $ with $\Phi$ smooth constant for large $x$. We will denote $\Delta_i$, $\nabla_i$ and $d\omega_i$ the Laplacian, gradient and area corresponding to this extension on $\R\times \Tu$. Since $u_i\in L^2(\Omega_i)$, we can consider the function $u_i$ in this local coordinates (in what follows we will not distinguish $u_i$ with its representent in local coordinates) and make the computation (in the sense of distributions of $\Omega$ and of $\R^*\times \Tu$ in local coordinates):
\bna
-\Delta_i u_i+u_i= -2\nabla_i u \cdot \nabla_i \Psi_i -(\Delta_i \Psi_i)u.
\ena
Remark that in these local coordinates $u$ has only a true meaning for small $x$ but we can then extend this equality on $\R^*\times \Tu$ since $\Psi_i$ is compactly supported. Moreover, $\nabla_i \Psi_i$ and $\Delta_i \Psi_i$ is supported outside of the zone of degeneracy $Z$ where $u$ is $C^{\infty}$. So, we get $-\Delta_i u_i+u_i \in C^{\infty}_0(\R^*\times \Tu)$ and so $-\Delta_i u_i+u_i\in L^2(M,d\omega_i)$. This is true in the sense of distributions, but that means that for any $\varphi \in C^{\infty}_0(\R^*\times \Tu)$, we have 
\bna
\left|(u_i,-\Delta_i \varphi +\varphi)_{L^2(\R^*\times \Tu,d\omega_i)}\right|\leq C \nor{\varphi}{L^2(\R^*\times \Tu,d\omega_i)}.
\ena   
In particular, that means that $u_i$ belongs to the domain of the adjoint of the operator $-\Delta_i+I$ with domain $C^{\infty}_0(\R^*\times \Tu)$. But Theorem \ref{Grushinvaress} gives that this operator is essentially self-adjoint so, $u_i \in D(L_i)$ where $L_i$ is the self-adjoint (Friedrichs) extension of $-\Delta_i+I$ on $\R^*\times \Tu$. In particular, we can write (the right hand side has to be understood as a limit for a sequence $u_{i,n}$ in $C^{\infty}_0(\Omega)$ converging strongly to $u_i$ for the norm of the quadratic form)  
\bna
(u_i,-\Delta_i u_i +u_i)_{L^2(\R^*\times \Tu,d\omega_i)}=\nor{\nabla_i u_i}{L^2(\R^*\times \Tu,d\omega_i)}^2+\nor{u_i}{L^2(\R^*\times \Tu,d\omega_i)}^2.
\ena
This quantity also makes sense when $u_i$ is considered on the manifold $M$ and the extension of the metric was chosen so that $-\Delta_i u_i +u_i=-\Delta u_i +u_i$ in the sense of ditributions and $\nabla_iu_i=\nabla u_i$.
 
The same result holds for $i\in I$ corresponding to a Riemannian zone.

Moreover, if $i\neq j \in I$, the common support of $u_i$ and $u_j$ do not intersect the degeneracy zone $Z$. So, these functions are $C^{\infty}$ in this zone and we can write
\bna
(u_i,-\Delta u_j +u_j)_{L^2(M)}=(\nabla u_i,\nabla u_j)_{L^2(M)}+(u_i,u_j)_{L^2(M)}
\ena
By summing up, we get 
\bna
0=(u,-\Delta u +u)_{L^2(M,d\omega)}=(\nabla u,\nabla u)_{L^2(M,d\omega)}+(u,u)_{L^2(M,d\omega)}
\ena
and $u\equiv 0$.
\enp
From now on, we will use the same notation $\Delta$ for the self-adjoint extension of the symmetric operator $\Delta$.

 To prove the second statement of Theorem \ref{t-1} we have to prove the following.
\bp
The domain of $\Delta$ is given by
\bnan
\label{domainDeltadistrib}
D(\Delta)=\left\lbrace u\in L^2(M,d\omega) \left| \Delta_{D,g} u \in L^2(M,d\omega)\right.\right\rbrace 
\enan
 where $\Delta_{D,g} u$ is $\Delta u$ seen as a distribution in $\Omega=M\setminus \Zz$. 
\ep
\bnp
Let $T$ be the Laplace operator defined with domain given by (\ref{domainDeltadistrib}). We easily see that $D(T)=\left\lbrace u\in L^2(M,d\omega) \left| \exists C :~\forall \varphi \in C^{\infty}_0(\Omega), |(u,\Delta \varphi)_{L^2}|\leq C \norL{\varphi} \right. \right\rbrace =D(\Delta^*)=D(\Delta)$.

\enp
\noindent
 To prove the third statement of Theorem \ref{t-1} we have to prove the following.
\bp
Denote by $\Delta$ the positive self-adjoint operator defined by  Proposition \ref{p-thmesscompact}, then the resolvent $(-\Delta +1)^{-1}$ is well defined and compact from $L^2(M,d\omega)$ to itself with the measure defined by the metric. 

Therefore, the spectrum of $\Delta$ is discrete and consists of eigenvalues of finite multiplicity.
\ep
\bnp
The fact that $(-\Delta +1)^{-1}$ is well defined comes from the positivity of $-\Delta$.  Now
let $u_n$ be a bounded sequence in $L^2(M,d\omega)$ and $f_n=(-\Delta +1)^{-1}u_n$. By density, we can assume $u_n$ in $C^{\infty}_0(M\setminus \Zz)$ and $f_n\in C^{\infty}(M\setminus \Zz)$ by elliptic regularity. Then, we have $\int_M |\nabla f_n|^2+|f_n|^2=(u_n,f_n)_{L^2} $ bounded. 

Consider the partition of unity $(\Psi_i)_{i\in I}$ introduced in the proof of the previous Theorem and denote $f_{i,n}:=\Psi_i f_n$. 

Let $i\in I$ be an index corresponding to a ``Grushin zone''.  By the formula 
\bna
(-\Delta+I)f_{n,i}=\Psi_i u_n-2\nabla \Psi_i \cdot \nabla f_n -(\Delta \Psi_i)f_n ,
\ena
we have that $(-\Delta+I)f_{n,i}$ is bounded in $L^2(M)$. By using Lemma \ref{lmcompactcoord}, we get that $\rho(x) (-\Delta_i+1)^{-1}(-\Delta+I)f_{n,i}$ is compact for a function $\rho\in C^{\infty}_0(\R)$ defined in some local coordinate charts. Here, we have denoted $\Delta_i$ the Laplacian for an extension of the local metric to $\R\times \Tu$. To finish, we only have to notice that $(-\Delta_i+1)^{-1}(-\Delta+I)f_{n,i}=(-\Delta_i+1)^{-1}(-\Delta_i+1)f_{n,i}=f_{n,i}$ because of the support of $f_{n,i}$ (note that it is not the case for $(-\Delta_i+1)^{-1}$ because the resolvent depends on the extension). 

So, choosing $\rho$ such that $\rho(x) f_{n,i}=f_{n,i}$ in local coordinates, we have proved that each sequence $f_{n,i}$ is compact for a Grushin zone. The same result holds for a Riemannian zone because the equivalent of Lemma \ref{lmcompactcoord} still holds for any Riemannian extension. By summing up, we get that the sequence $f_n$ is convergent, up to extraction, which yields the result.
\enp
\appendix
\section{The Martinet case}
In this Section, we prove Theorem \ref{t-martinet}. 
\ugo{A sub-Riemannian structure is a triple $(M,\bD,{\bf g})$, where $M$ is  a smooth manifold, $\bD$ is a smooth vector distribution of constant rank satisfying the H\"ormander condition, and ${\bf g}$ is a Riemannian metric on $\bD$. 
Let $\bD_1 = \bD$ and $\bD_{k+1}=\bD_k+[\bD,\bD_k]$. A sub-Riemannian structure is said to be equiregular if the dimension of $\bD_k$ does not depend on the point.

Almost-Riemannian and sub-Riemannian structures can be treated in the unified setting of rank-varying sub-Riemannian structures, see \cite[Chapter 3]{AgrBarBoscbook} and \cite[Definition 2]{euler}. Here for sake of readability we omit this point of view.

In this section, we briefly treat the Martinet sub-Riemannian manifold (see for instance \cite{bonnard-book,Montgomery}), defined by $M=\R^3$, $\bD(q)=span\{X_1(q),X_2(q)\}$, and ${\bf g}(X_i,X_j)=\delta_{ij}$, $i,j=1,2$, }
where
\bqn
X_1=
\left(\ba{c}
1\\0\\\frac{y^2}{2}
\ea
\right),~~~X_2=
\left(\ba{c}
0\\1\\0
\ea
\right).
\eqn
Notice that $X_3=[X_1,X_2]=\left(\ba{c}
0\\0\\-y
\ea
\right)$.
\ugo{Hence $X_1$, $X_2$ and $X_3$ span the tangent space at any point outside the so called Martinet plane  
$\Zz=\left\lbrace y=0\right\rbrace $, which is the region where the structure is not equiregular.} However $X_4=[[X_1,X_2],X_2]=[X_3,X_2]=\left(\ba{c}
0\\0\\1
\ea
\right)$, hence the H\"ormander condition is fulfilled on the whole space. 

With this structure, $\R^3\setminus\Zz$ is a $3D$ sub-Riemannian contact manifold.  On such a structure it is possible to define intrinsically a volume form which is given by $dX_1 \wedge dX_2 \wedge dX_3$, where $\{dX_1, dX_2, dX_3\}$ is the dual basis to $\{X_1,X_2,X_3\}$. Such a volume form is independent on the choice of the orthonormal basis which define the sub-Riemannian structure, and it is called the Popp measure (see  \cite{Montgomery} and Proposition 8 of \cite{HeatKernel}).
One gets $dX_1=dx$, $dX_2=dy$ and $dX_3=\frac{y}{2}dx-\frac{1}{y}dz$ so that $dX_1 \wedge dX_2 \wedge dX_3=- \frac{1}{y}dx\wedge dy\wedge dz$ and the corresponding density is $d\omega=\frac{1}{|y|}dxdydz$. This allows to define the sub-Riemannian Laplacian as the divergence of the sub-Riemannian gradient defined as in formula \r{ars-grad} (see Remark 14 of \cite{HeatKernel}),
\bqn
\Delta_{sr}=(X_1)^2+(X_2)^2-\frac{1}{y}X_2=(\partial_x+\frac{y^2}{2}\partial_z)^2+\partial_y^2-\frac{1}{y}\partial_y.
\eqn
 
Notice that the singularity in the first order term of  $\Delta_{sr}$ appears similarly as it does in the Grushin case.
We do the same reasoning.

For simplicity, we compactify in $x$ and $z$ and consider the same structure on $M=\Tu_x\times \R_y \times \Tu_z$. We will prove that $\Delta_{sr}$ is essentially self-adjoint on $L^2(M,d\omega)$ with domain $C^{\infty}_0(M\setminus \Zz)$. Again, we make the change of variable $f=\sqrt{|y|}g$ so that 
\bqn
\Delta_{sr}f=\sqrt{|y|}\left( (\partial_x+\frac{y^2}{2}\partial_z)^2g+\partial_y^2g -\frac{3}{4|y|^2}g\right)=\sqrt{|y|}\widetilde{\Delta_{sr}}g.
\eqn
So, we are left to prove that $\widetilde{\Delta_{sr}}=(\partial_x+\frac{y^2}{2}\partial_z)^2+\partial_y^2 -\frac{3}{4y^2}$ is essentially self-adjoint on $L^2(\Tu\times\R\times \Tu,dx\,dy\,dz)$ with domain $C^{\infty}_0(M\setminus \Zz)$.
We compute the Fourier transform in $x$ and $z$ and get the decomposition 
\bna
L^2(\Tu\times \R\times \Tu)\approx \oplus_{(k,l)\in \Z^2}^{\perp} H_{(k,l)}
\ena
where $H_{(k,l)}=L^2(\R)$ and the operator $L$ acts on each $H_{(k,l)}$ by 
\bna
(L u)_{(k,l)}(y)=\partial_y^2 u_{(k,l)}-(k+\frac{y^2}{2}l)^2 u_{(k,l)} -\frac{3}{4y^2} u_{(k,l)}=\partial_y^2 u_{(k,l)}-V_{(k,l)} u_{(k,l)}.
\ena
with $ V_{(k,l)}(y)\geq \frac{3}{4y^2}$. We conclude as in the Grushin case.

\ugo{This result suggests the general conjecture that for  a sub-Riemannian structure which is rank-varying or not equiregular on an hypersurface, the singular set acts as a barrier for the heat flow and for a quantum particle.}

\bibliographystyle{abbrv}
\bibliography{biblio-grushin-beltrami-v-2.7}
}

\end{document}

%% file: fig-intro-ugo.pstex_t
\begin{picture}(0,0)%
\includegraphics{fig-intro-ugo.pstex}%
\end{picture}%
\setlength{\unitlength}{1776sp}%
\begingroup\makeatletter\ifx\SetFigFont\undefined%
\gdef\SetFigFont#1#2#3#4#5{%
  \reset@font\fontsize{#1}{#2pt}%
  \fontfamily{#3}\fontseries{#4}\fontshape{#5}%
  \selectfont}%
\fi\endgroup%
\begin{picture}(12999,9624)(3589,-10273)
\put(5926,-7486){\makebox(0,0)[lb]{\smash{{\SetFigFont{9}{10.8}{\rmdefault}{\mddefault}{\updefault}{\color[rgb]{0,0,0}(F3): Tangency points}%
}}}}
\put(4126,-2536){\makebox(0,0)[lb]{\smash{{\SetFigFont{9}{10.8}{\rmdefault}{\mddefault}{\updefault}$X_1=\frp{}{x}\quad X_2=e^{\phi(x,y)}\frp{}{y}$}}}}
\put(5926,-1636){\makebox(0,0)[lb]{\smash{{\SetFigFont{9}{10.8}{\rmdefault}{\mddefault}{\updefault}{\color[rgb]{0,0,0}(F1): Riemannian points}%
}}}}
\put(6001,-4486){\makebox(0,0)[lb]{\smash{{\SetFigFont{9}{10.8}{\rmdefault}{\mddefault}{\updefault}{\color[rgb]{0,0,0}(F2): Grushin points}%
}}}}
\put(4951,-9736){\makebox(0,0)[lb]{\smash{{\SetFigFont{9}{10.8}{\rmdefault}{\mddefault}{\updefault}{\color[rgb]{0,0,0}$\psi,\xi\in\con^{\infty},\,\,\,\psi(0)\neq0$ }%
}}}}
\put(3976,-9061){\makebox(0,0)[lb]{\smash{{\SetFigFont{9}{10.8}{\rmdefault}{\mddefault}{\updefault}$X_1=\frp{}{x}\quad X_2=(y-x^2\psi(x))e^{\xi(x,y)}\frp{}{y}$}}}}
\put(15076,-9811){\makebox(0,0)[lb]{\smash{{\SetFigFont{9}{10.8}{\rmdefault}{\mddefault}{\updefault}{\color[rgb]{1,0,0}$X_2$}%
}}}}
\put(15451,-9436){\makebox(0,0)[lb]{\smash{{\SetFigFont{9}{10.8}{\rmdefault}{\mddefault}{\updefault}{\color[rgb]{0,0,1}$X_1$}%
}}}}
\put(12226,-7786){\makebox(0,0)[lb]{\smash{{\SetFigFont{9}{10.8}{\rmdefault}{\mddefault}{\updefault}{\color[rgb]{0,0,0}$\Zz$}%
}}}}
\put(11401,-4336){\makebox(0,0)[lb]{\smash{{\SetFigFont{9}{10.8}{\rmdefault}{\mddefault}{\updefault}{\color[rgb]{0,0,0}$\Zz$}%
}}}}
\put(13351,-4411){\makebox(0,0)[lb]{\smash{{\SetFigFont{9}{10.8}{\rmdefault}{\mddefault}{\updefault}{\color[rgb]{1,0,0}$X_2$}%
}}}}
\put(13426,-5011){\makebox(0,0)[lb]{\smash{{\SetFigFont{9}{10.8}{\rmdefault}{\mddefault}{\updefault}{\color[rgb]{0,0,1}$X_1$}%
}}}}
\put(3976,-5536){\makebox(0,0)[lb]{\smash{{\SetFigFont{9}{10.8}{\rmdefault}{\mddefault}{\updefault}$X_1=\frp{}{x}\quad X_2=xe^{\phi(x,y)}\frp{}{y}$}}}}
\put(4501,-6361){\makebox(0,0)[lb]{\smash{{\SetFigFont{9}{10.8}{\rmdefault}{\mddefault}{\updefault}{\color[rgb]{0,0,0}$\phi\in\con^{\infty},\,\,\phi(0,y)\equiv 0$ }%
}}}}
\put(13426,-2086){\makebox(0,0)[lb]{\smash{{\SetFigFont{9}{10.8}{\rmdefault}{\mddefault}{\updefault}{\color[rgb]{0,0,1}$X_1$}%
}}}}
\put(13351,-1486){\makebox(0,0)[lb]{\smash{{\SetFigFont{9}{10.8}{\rmdefault}{\mddefault}{\updefault}{\color[rgb]{1,0,0}$X_2$}%
}}}}
\put(4576,-3436){\makebox(0,0)[lb]{\smash{{\SetFigFont{9}{10.8}{\rmdefault}{\mddefault}{\updefault}{\color[rgb]{0,0,0}$\phi\in\con^{\infty},\,\,\phi(0,y)\equiv 0$ }%
}}}}
\end{picture}%